\documentclass[journal]{IEEEtran}
\usepackage{cite}

\usepackage[pdftex]{graphicx}
\usepackage{tikz,pgfplots}
\usetikzlibrary{positioning}
\pgfplotsset{compat=1.14} 
\tikzstyle{vertex} = [fill,shape=circle,node distance=30pt]
\tikzstyle{edge} = [fill,opacity=.6,fill opacity=.5,line cap=round, line join=round, line width=10pt]
\tikzstyle{elabel} =  [fill,shape=circle,node distance=30pt,fill opacity=.9]
\definecolor{mygray}{gray}{0.95}
\pgfdeclarelayer{background}
\pgfsetlayers{background,main}
\usepackage[most]{tcolorbox}
\definecolor{mypurple}{rgb}{0.59, 0.44, 0.84}
\ifCLASSINFOpdf
\else
\fi
\usepackage{amsmath}
\usepackage{amssymb}
\usepackage{mathrsfs} 
\usepackage{mathtools}
\usepackage[utf8]{inputenc}
\usepackage[english]{babel}
\usepackage{amsthm}

\usepackage{algorithmic}

\usepackage[ruled]{algorithm}

\usepackage{array}
\begin{document}
%
\title{Controllability of Hypergraphs}
%
%
%

\author{Can~Chen,
        Amit~Surana,~\IEEEmembership{Member,~IEEE,} Anthony~M.~Bloch,~\IEEEmembership{Fellow,~IEEE,} 
        and~Indika~Rajapakse
\thanks{Manuscript received October 23, 2020; revised February 18, 2021. This work is supported in part by AFOSR Award No: FA9550-18-1-0028, in part by NSF grant DMS 1613819, in part by the Smale Institute, and in part by the Guaranteeing AI Robustness Against Deception (GARD) program from DARPA (I2O).} 
\thanks{C. Chen is with the Department of Mathematics and the Department of Electrical Engineering and Computer Science, University of Michigan, Ann Arbor, MI 48109 USA (e-mail: canc@umich.edu).}
\thanks{A. Surana is with Raytheon Technologies Research Center, East Hartford, CT 06108 USA (e-mail: amit.surana@rtx.com).}
\thanks{A. Bloch is with the Department of Mathematics, University of Michigan, Ann Arbor, MI 48109 USA (e-mail: abloch@umich.edu).}
\thanks{I. Rajapakse is with the Department of Computational Medicine \& Bioinformatics, Medical School and the Department of Mathematics, University of Michigan, Ann Arbor, MI 48109 USA (e-mail: indikar@umich.edu).}
}

%
%

\markboth{IEEE TRANSACTIONS ON NETWORK SCIENCE AND ENGINEERING}%
{Chen \MakeLowercase{\textit{et al.}}}
%



\maketitle

\begin{abstract}

In this paper, we develop a notion of controllability for hypergraphs via tensor algebra and polynomial control theory. Inspired by uniform hypergraphs, we propose a new tensor-based multilinear dynamical system representation, and derive a Kalman-rank-like condition to determine the minimum number of control nodes (MCN) needed to achieve controllability of even uniform hypergraphs. We present an efficient heuristic to obtain the MCN. MCN can be used as a measure of robustness, and we show that it is related to the hypergraph degree distribution in simulated examples.  Finally, we use MCN to examine robustness in real biological networks.

\end{abstract}


%
\IEEEpeerreviewmaketitle

\section{Introduction}
\IEEEPARstart{M}{any} complex systems are studied using a network perspective, which offers unique insights in social sciences, cell biology, neuroscience and computer science \cite{strogatz2001exploring,Amaral11149,scale_free, small_world2,small_world}. For example, recent advances in genomics technology, such as genome-wide chromosomal conformation capture (Hi-C), have inspired us to consider the human genome as a dynamic network \cite{Rajapakse711,RIED20171,9119161}.  Studying such dynamic networks often requires introducing external inputs into the networks in order to steer the network dynamics towards a desired state \cite{1100557,1023188,tanner2004controllability,Rajapakse17257}.  This process is akin to the notion of controllability in classical control theory. A dynamical system is controllable if it can be driven from any initial state to any target state within a finite time given a suitable choice of control inputs. 

Controlling complex networks is one of the most challenging  problems in modern network science \cite{doi.org/10.1038/nature10011,doi.org/10.1038/ncomms3447,gates2016control,zanudo2017structure,8851061,wang2020controllability,wang2017physical}. Lin \cite{1100557} first proposed the concept of structural controllability of directed graphs in 1970s. Later on, Tanner \cite{tanner2004controllability} and Rahmani \textit{et al.} \cite{rahmani2009controllability,rahmani2006controlled} applied the idea of structural controllability for multi-agent systems with the aim of selecting a subset of agents (called leaders) which are able to control the whole system by exploiting the graph Laplacian and linear control theory. In particular, Rahmani \textit{et al.} \cite{rahmani2009controllability} also showed how the symmetry structure of a graph directly relates to the controllability of the corresponding multi-agent system. 

In 2011, Liu \textit{et al.} \cite{doi.org/10.1038/nature10011,RevModPhys.88.035006} explored the (structural) controllability of complex graphs with $n$ nodes by using the canonical linear time-invariant dynamics
\begin{equation}\label{eq:1}
\dot{\textbf{x}} = \textbf{A}\textbf{x} + \textbf{B}\textbf{u},
\end{equation}
where, $\textbf{A}\in\mathbb{R}^{n\times n}$ is the adjacency matrix of a graph, and $\textbf{B}\in\mathbb{R}^{n\times m}$ is the control matrix.  The time-dependent vector $\textbf{x}\in\mathbb{R}^n$ captures the states of the nodes, and $\textbf{u}\in\mathbb{R}^m$ is the time-dependent control vector. The authors exploited the Kalman rank condition, i.e., the linear system (\ref{eq:1}) is controllable if and only if the controllability matrix
\begin{equation}
\textbf{C} = \begin{bmatrix} \textbf{B} & \textbf{A}\textbf{B} & \dots & \textbf{A}^{n-1}\textbf{B}\end{bmatrix}
\end{equation}
has full rank, to determine the minimum number of control nodes (MCN) in order to achieve full control of the graph (similar to the role of leaders discussed previously). They also identified the MCN of a graph using the idea of ``maximum matching''  \cite{doi.org/10.1038/nature10011}. In addition, Yuan \textit{et al.} \cite{doi.org/10.1038/ncomms3447} developed a notion of exact controllability of graphs. They took advantage of the Popov-Belevitch-Hautus rank condition (i.e., the linear system (\ref{eq:1}) is controllable if and only if $\text{rank}(\begin{bmatrix}s\textbf{I}-\textbf{A} & \textbf{B}\end{bmatrix})=n$ for any complex number $s$) to prove that for an arbitrary graph, the MCN is determined by the maximum geometric multiplicity of the eigenvalues of the corresponding adjacency matrix \textbf{A}. Furthermore, Nacher \textit{et al.} \cite{nacher2019finding} analyzed MCN required to fully control multilayer graphs, and a similar notion of exact controllability for multilayer graphs is defined in \cite{yuan2014exact}. 

However, most real world data representations are multidimensional, and using graph (or even multilayer graph) models to describe them may result in a loss of higher-order information \cite{9119161,7761624}. A hypergraph is a generalization of a graph in which its hyperedges can join any number of nodes \cite{hypergraph1,chung1982minimal,8789484,9222341}. Thus, hypergraphs can capture multidimensional relationships unambiguously \cite{7761624}. Examples of hypergraphs include co-authorship networks, film actor/actress networks, and protein-protein interaction networks \cite{newman2010networks}. Moreover, a hypergraph can be represented by a tensor if its hyperedges contain the same number of nodes, referred to as a uniform hypergraph \cite{9119161}. Tensors are multidimensional arrays generalized from vectors and matrices that preserve multidimensional patterns and capture higher-order interactions and coupling within multiway data \cite{chen_2020}. The dynamics of uniform hypergraphs thus can be naturally described by a tensor-based multilinear system (multilinear in the sense of multilinear algebra).

Multilinear dynamical systems were first proposed by Rogers \textit{et al.} \cite{rogers_2013} and Surana \textit{et al.} \cite{7798500} where the system evolution is generated by the action of multilinear operators which are formed using Tucker products of matrices. Chen \textit{et al.} \cite{doi:10.1137/1.9781611975758.18} developed the tensor algebraic conditions for stability, reachability and observability for input/output multilinear time-invariant systems. By using tensor unfolding, an operation that transforms a tensor into a matrix, a multilinear system can be unfolded to a corresponding linear system. However, the tensor-based multilinear systems, proposed in this paper, are different from the ones defined in \cite{rogers_2013, 7798500, doi:10.1137/1.9781611975758.18}, and in fact they belong to the family of nonlinear polynomial systems. Hence, they can capture network dynamics more precisely than systems based on standard graphs which use linear dynamics assumption. Basic knowledge of nonlinear control concepts such as Lie algebra and Lie brackets is required in order to better understand the controllability of such systems. The key contributions of this paper are as follows:
\begin{itemize}
\item We propose a new tensor-based multilinear system representation inspired by uniform hypergraphs, and study the controllability of such systems by exploiting tensor algebra and polynomial control theory. We establish a Kalman-rank-like condition to determine the controllability of even uniform hypergraphs. 
\item We establish theoretical results on the MCN of even uniform hyperchains, hyperrings and hyperstars as well as complete even uniform hypergraphs. We also observe that the MCN of odd uniform hypergraphs, identified by the Kalman-rank-like condition, behaves similarly to that of even uniform hypergraphs in simulated examples (although the condition is not applicable in terms of controllability). We discover that MCN is related to the hypergraph degree distribution, and high degree nodes are preferred to be controlled in these configurations and their variants. 

\item We propose MCN as a measure of robustness for uniform hypergraphs, and use it to quantify structural differences. We present applications to two real world biological examples: a  mouse  neuron  endomicroscopy dataset and an allele-specific Hi-C dataset. 
\item We present a fast and memory-efficient computational framework for determining the rank of a matrix related to the Kalman-rank-like condition for nonlinear controllability. In addition, we propose a heuristic approach to identify the MCN of uniform hypergraphs efficiently.
\item We perform preliminary explorations of the controllability of general non-uniform hypergraphs.
\end{itemize}

The paper is organized into seven sections. We start with the basics of tensor algebra in Section \ref{sec:2.1}. In Section \ref{sec:2.2}, we discuss the notion of uniform hypergraphs and extend the definitions of chains, rings and stars from graph theory to uniform hypergraphs. We propose a new tensor-based multilinear system to capture the dynamics of uniform hypergraphs in Section \ref{sec:2.3}. We then formulate a Kalman-rank-like condition to determine the controllability of even uniform hypergraphs in Section \ref{sec:2.4}. We establish theoretical results on the MCN of even uniform hyperchains, hyperrings and hyperstars as well as complete even uniform hypergraphs in Section \ref{sec:2.6}. In Section \ref{sec:2.7}, we argue that MCN can be used to measure hypergraph robustness, and provide a heuristic approach to find the MCN efficiently. Five numerical examples are presented in Section \ref{sec:3}.  Finally, we discuss the controllability of general hypergraphs in Section \ref{sec:4}  and conclude with future directions in Section \ref{sec:5}. 

\section{Preliminaries}
\subsection{Tensors}\label{sec:2.1}
We take most of the concepts and notations for tensor algebra from the comprehensive works of Kolda \textit{et al.} \cite{doi:10.1137/07070111X, Kolda06multilinearoperators}. A \textit{tensor} is a multidimensional array. The \textit{order} of a tensor is the number of its dimensions,  and each dimension is called a \textit{mode}. A $k$-th order tensor usually is denoted by $\textsf{T}\in \mathbb{R}^{n_1\times n_2\times  \dots \times n_k}$.  It is therefore reasonable to consider scalars $x\in\mathbb{R}$ as zero-order tensors, vectors $\textbf{v}\in\mathbb{R}^{n}$ as first-order tensors, and matrices $\textbf{M}\in\mathbb{R}^{m\times n}$ as second-order tensors. A tensor is called \textit{cubical} if every mode is the same size, i.e., $\textsf{T}\in \mathbb{R}^{n\times n\times  \dots \times n}$. A cubical tensor \textsf{T} is called \textit{supersymmetric} if $\textsf{T}_{j_1j_2\dots j_k}$ is invariant under any permutation of the indices.

The \textit{tensor vector multiplication} $\textsf{T} \times_{p} \textbf{v}$ along mode $p$ for a vector $\textbf{v}\in  \mathbb{R}^{n_p}$ is defined  by
\begin{equation}
(\textsf{T} \times_{p} \textbf{v})_{j_1j_2\dots j_{p-1}j_{p+1}\dots j_k}=\sum_{j_p=1}^{n_p}\textsf{T}_{j_1j_2\dots j_p\dots j_k}\textbf{v}_{j_p}, 
\end{equation}
which can be extended to 
\begin{equation}\label{eq6}
\begin{split}
\textsf{T}\times_1 \textbf{v}_1 \times_2\textbf{v}_2\times_3\textbf{v}_3\dots \times_{k}\textbf{v}_k=\textsf{T}\textbf{v}_1\textbf{v}_2\textbf{v}_3\dots\textbf{v}_k\in\mathbb{R}
\end{split}
\end{equation}
for $\textbf{v}_p\in \mathbb{R}^{ n_p}$. The expression (\ref{eq6}) is also known as the homogeneous polynomial associated with \textsf{T}. If $\textbf{v}_p=\textbf{v}$ for all $p$, we write  (\ref{eq6})  as $\textsf{T}\textbf{v}^{k}$ for simplicity. 

\subsection{Uniform Hypergraphs}\label{sec:2.2}
We borrow some fundamental concepts of hypergraphs from the works \cite{BANERJEE201714,COOPER20123268,HU2014140,qi2013}. An \textit{undirected hypergraph} \textsf{G} is a pair such that $\textsf{G} = \{\textbf{V}, \textbf{E}\}$ where $\textbf{V} = \{1,2,\dots,n\}$ is the node set and $\textbf{E} = \{e_1,e_2,\dots,e_p\}$ is the \textit{hyperedge} set with $e_l\subseteq \textbf{V}$ for $l=1,2,\dots,p$. Two nodes are called \textit{adjacent} if they are in the same hyperedge. A hypergraph is called \textit{connected} if, given two nodes, there is a path connecting them through hyperedges. If all hyperedges contain the same number of nodes, i.e., $|e_p| = k$ for $k\leq n$, \textsf{G} is called a \textit{$k$-uniform hypergraph}. Here $|\cdot|$ denotes the cardinality of a set. A $k$-uniform hypergraph can be represented by a $k$-th order $n$-dimensional supersymmetric tensor.

\textit{Definition 1 (\cite{BANERJEE201714,COOPER20123268,HU2014140,qi2013}):} Let \textsf{G} = \{\textbf{V}, \textbf{E}\} be a $k$-uniform hypergraph with $n$ nodes. The \textit{adjacency tensor} $\textsf{A}\in\mathbb{R}^{n\times n\times \dots\times n}$ of \textsf{G}, which is a $k$-th order $n$-dimensional supersymmetric tensor, is defined as
\begin{equation}\label{eq:98}
\textsf{A}_{j_1j_2\dots j_k} = \begin{cases} \frac{1}{(k-1)!} \text{ if $(j_1,j_2,\dots,j_k)\in \textbf{E}$}\\ \\0, \text{ otherwise}\end{cases}.
\end{equation}

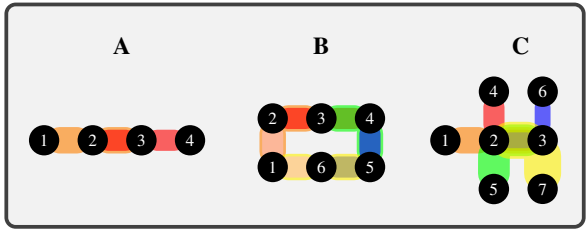
\begin{figure}[t]
\centering
\tcbox[colback=mygray,top=3pt,left=5pt,right=5pt,bottom=5pt]{
\begin{tikzpicture}[scale=0.92, transform shape]
\node[vertex,text=white,scale=0.7] (v1) {1};
\node[vertex,right of=v1,node distance = 20pt,text=white,scale=0.7] (v2) {2};
\node[vertex,right of=v2,node distance = 20pt,text=white,scale=0.7] (v3) {3};
\node[vertex,right of=v3,node distance = 20pt,text=white,scale=0.7] (v4) {4};

\node[vertex,below right = 2pt and 25pt of v4,node distance = 30pt,text=white,scale=0.7] (v10) {1};
\node[vertex,above of=v10,node distance = 20pt,text=white,scale=0.7] (v11) {2};
\node[vertex,right of=v11,node distance = 20pt,text=white,scale=0.7] (v12) {3};
\node[vertex,right of=v12,node distance = 20pt,text=white,scale=0.7] (v13) {4};
\node[vertex,below of=v13,node distance = 20pt,text=white,scale=0.7] (v14) {5};
\node[vertex,left of=v14,node distance = 20pt,text=white,scale=0.7] (v15) {6};

\node[vertex,above of=v12,mygray, text=black,scale=1] (b) {\textbf{B}};
\node[vertex,left of=b,mygray, text=black,scale=1,node distance = 82pt] (a) {\textbf{A}};
\node[vertex,right of=b,mygray, text=black,scale=1,node distance = 82pt] (c) {\textbf{C}};

\node[vertex,right of=v4,text=white,scale=0.7,node distance=105pt] (v22) {1};
\node[vertex,right of=v22,text=white,scale=0.7,node distance = 20pt] (v23) {2};
\node[vertex,right of=v23,text=white,scale=0.7,node distance = 20pt] (v24) {3};
\node[vertex,above left = 11pt and 11pt of v24,text=white,scale=0.7] (v25) {4};
\node[vertex,below left = 11pt and 11pt of v24,text=white,scale=0.7] (v26) {5};
\node[vertex,above of=v24,text=white,scale=0.7,node distance = 20pt] (v27) {6};
\node[vertex,below of=v24,text=white,scale=0.7,node distance = 20pt] (v28) {7};

\begin{pgfonlayer}{background}
\draw[edge,color=orange] (v1) -- (v2) -- (v3);
\draw[edge,color=red, line width = 8] (v2) -- (v3) -- (v4);

\draw[edge,color=orange] (v10) -- (v11) -- (v12);
\draw[edge,color=red, line width = 8] (v11) -- (v12) -- (v13);
\draw[edge,color=green] (v12) -- (v13) -- (v14);
\draw[edge,color=blue, line width=8] (v13) -- (v14) -- (v15);
\draw[edge,color=yellow] (v14) -- (v15) -- (v10);
\draw[edge,color=pink,line width=8] (v15) -- (v10) -- (v11);

\draw[edge,color=orange] (v22) -- (v23) -- (v24);
\draw[edge,color=red, line width=8] (v25) -- (v23) -- (v24);
\draw[edge,color=green, line width=12] (v26) -- (v23) -- (v24);
\draw[edge,color=blue, line width=6] (v23) -- (v24) -- (v27);
\draw[edge,color=yellow, line width=14] (v23) -- (v24) -- (v28);

\end{pgfonlayer}
\end{tikzpicture}
}
\caption{Examples of hyperchains, hyperrings and hyperstars. (A) 3-uniform hyperchain with $e_1=\{1,2,3\}$ and $e_2=\{2,3,4\}$. (B) 3-uniform hyperring with $e_1=\{1,2,3\}$, $e_2=\{2,3,4\}$, $e_3=\{3,4,5\}$, $e_4=\{4,5,6\}$, $e_5=\{5, 6, 1\}$ and $e_6 =\{6,1,2\}$. (C) 3-uniform hyperstar with $e_1=\{1,2,3\}$, $e_2=\{2,3,4\}$, $e_3=\{2,3,5\}$, $e_4=\{2,3,6\}$ and $e_5=\{2,3,7\}$.}
\label{fig:10}
\end{figure}

Similarly to standard graphs, the \textit{degree} of node $j$ of a uniform hypergraph is defined as 
\begin{equation}\label{eq:r6}
    d_j=\sum_{j_2=1}^n\sum_{j_3=1}^n\dots \sum_{j_k=1}^n\textsf{A}_{jj_2j_3\dots j_k}.
\end{equation}
Note that the choice of the nonzero coefficient $\frac{1}{(k-1)!}$ in (\ref{eq:98}) guarantees that the degree of each node is equal to the number of hyperedges that contain that node, which is consistent with the notion of degree in standard graphs. Moreover, the controllability framework  we develop in Section \ref{sec:2.4} can be applied to arbitrary weighted uniform hypergraphs, see Section \ref{sec:4}. The \textit{degree distribution} of a hypergraph is the distribution of the degrees over all nodes. If all nodes have the same degree $d$, then \textsf{G} is called \textit{$d$-regular}. Given any $k$ nodes, if they are contained in one hyperedge, then \textsf{G} is called \textit{complete}. In the following, we extend the definitions of chains, rings and stars from graph theory to uniform hypergraphs.

\textit{Definition 2:} A $k$-uniform \textit{hyperchain} is a sequence of $n$ nodes such that every $k$ consecutive nodes are adjacent, i.e., nodes $j,j+1,\dots,j+k-1$ are contained in one hyperedge for $j=1,2,\dots,n-k+1$. 

\textit{Definition 3:} A $k$-uniform \textit{hyperring} is a sequence of $n$ nodes such that every $k$ consecutive nodes are adjacent, i.e., nodes $\sigma_n(j),\sigma_n(j+1),\dots,\sigma_n(j+k-1)$ are contained in one hyperedge for $j=1,2,\dots,n$, where $\sigma_n(j)=j$ for $j\leq n$ and $\sigma_n(j)=j-n$ for $j>n$.

\textit{Definition 4:} A $k$-uniform \textit{hyperstar} is a collection of $k-1$ internal nodes that are contained in all the hyperedges, and $n-k+1$ leaf nodes such that every leaf node is contained in one hyperedge with the internal nodes.

In $k$-uniform hyperchains, hyperrings and hyperstars, every two hyperedges have exactly $k-1$ overlapping nodes, see Fig. \ref{fig:10}. When $k=2$, they are reduced to standard chains, rings and stars.  We will determine the minimum number of control nodes (MCN) of uniform hyperchains, hyperrings and hyperstars in Section \ref{sec:2.6}.

\section{Hypergraph Controllability}\label{sec:2}
\subsection{Uniform Hypergraph Dynamics}\label{sec:2.3}

We represent the dynamics of a $k$-uniform hypergraph \textsf{G} with $n$ nodes by multilinear time-invariant differential equations using the adjacency tensor of \textsf{G}.

\textit{Definition 5:} Given a $k$-uniform hypergraph \textsf{G} with $n$ nodes, the dynamics of \textsf{G} with control inputs can be represented by 
\begin{equation}\label{eq:2}
\dot{\textbf{x}}= \textsf{A}\textbf{x}^{k-1} + \sum_{j=1}^m\textbf{b}_ju_j,
\end{equation}
where, $\textsf{A}\in\mathbb{R}^{n\times n\times \dots \times n}$ is the adjacency tensor of \textsf{G}, and $\textbf{B}=\begin{bmatrix}\textbf{b}_1 & \textbf{b}_2 & \dots & \textbf{b}_m\end{bmatrix}\in\mathbb{R}^{n\times m}$ is the control matrix. 

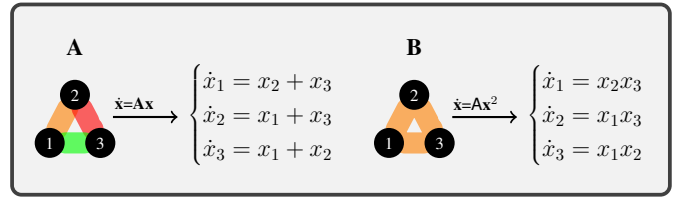
\begin{figure}[t!]
\centering
\tcbox[colback=mygray,top=5pt,left=5pt,right=-5pt,bottom=5pt]{
\begin{tikzpicture}[scale=0.91, transform shape]
\node[vertex,text=white,scale=0.7] (v1) {1};
\node[vertex,above right = 11pt and 1.5pt of v1,text=white,scale=0.7] (v2) {2};
\node[vertex,below right = 11pt and 1.5pt of v2,text=white,scale=0.7] (v3) {3};

\node[vertex,below right = 3pt and 40pt of v2,mygray,scale=0.5] (v7){};
\node[vertex,below right = 3pt and 6pt of v2,mygray,scale=0.5] (v8){};
\node[vertex,below right = -9pt and 12pt of v2,mygray,scale=0.7,text=black] (v9){  $\dot{\textbf{x}}$=$\textbf{A}\textbf{x}$};
\node[rectangle,below right = -20pt and 40pt of v2,mygray,scale=0.7,text=black] (v6) {
$\Large
\begin{cases}
\dot{x}_1 = x_2+x_3\\
\dot{x}_2 = x_1 + x_3\\
\dot{x}_3 = x_1 + x_2
\end{cases}
$};
\node[vertex,above of=v2,node distance=20pt,mygray,text=black,scale=0.7] (text1) {\Large\textbf{A}};
\path [->,shorten >=1pt,shorten <=1pt, thick](v8) edge node[left] {} (v7);

\node[vertex,right of=v3,node distance=120pt,text=white,scale=0.7](w1) {1};
\node[vertex,above right = 11pt and 1.5pt of w1,text=white,scale=0.7] (w2) {2};
\node[vertex,below right = 11pt and 1.5pt of w2,text=white,scale=0.7] (w3) {3};
\node[vertex,below right = 3pt and 40pt of w2,mygray,scale=0.5] (w7){};
\node[vertex,below right = 3pt and 6pt of w2,mygray,scale=0.5] (w8){};
\node[vertex,below right = -10pt and 12pt of w2,mygray,scale=0.7,text=black] (w9){  $\dot{\textbf{x}}$=$\textsf{A}\textbf{x}^{2}$};

\node[rectangle,below right = -20pt and 40pt of w2,mygray,scale=0.7,text=black] (w6) {
$\Large
\begin{cases}
\dot{x}_1 = x_2x_3\\
\dot{x}_2 = x_1x_3\\
\dot{x}_3 = x_1x_2
\end{cases}
$};
\node[vertex,above of=w2,node distance=20pt,mygray,text=black,scale=0.7] (text1) {\Large\textbf{B}};

\path [->,shorten >=1pt,shorten <=1pt, thick](w8) edge node[left] {} (w7);

\begin{pgfonlayer}{background}
\draw[edge,color=orange,line width=8pt] (v1) -- (v2);
\draw[edge,color=red,line width=8pt] (v2) -- (v3);
\draw[edge,color=green,line width=8pt] (v3) -- (v1);

\draw[edge,color=orange,line width=8pt] (w1) -- (w2) -- (w3) -- (w1);
\end{pgfonlayer}

\end{tikzpicture}
}
\caption{Graphs versus uniform hypergraphs. (A) Standard graph with three nodes and edges $e_1=\{1,2\}$, $e_2=\{2,3\}$ and $e_3=\{1,3\}$, and its corresponding linear dynamics. (B) 3-uniform hypergraph with three nodes and a hyperedge $e_1=\{1,2,3\}$, and its corresponding nonlinear dynamics.}
\label{fig:3}
\end{figure}

In this paper, we consider the case in which each input can only be imposed at one node, i.e., $\textbf{b}_j$ are the scaled standard basis vectors, similar to the treatments in \cite{doi.org/10.1038/nature10011,doi.org/10.1038/ncomms3447}. The time-dependent vector $\textbf{x}$ captures the state of the $n$ nodes, and the system is controlled using the time-dependent input vector $\textbf{u} = \begin{bmatrix} u_1 & u_2 & \dots & u_m\end{bmatrix}^\top\in\mathbb{R}^m$. The multilinear system (\ref{eq:2}) formulated by the tensor vector multiplications is indeed able to capture the simultaneous interactions among nodes for uniform hypergraphs as illustrated in Fig. \ref{fig:3}. All the interactions are characterized using multiplications instead of the additions that are typically used in a standard graph. It is known that  multiplication often stands for simultaneity and addition for sequentiality in many mathematical fields. For example, the probability of two independent events to happen at the same time is equal to the product of their individual probabilities. A similar form of nonlinear dynamical system representation has been used to model protein-protein interaction  \cite{fujarewicz2005fitting}, which as mentioned can be represented by hypergraphs. Compared to the dynamical systems of hypergraphs defined in \cite{timoteo2020,de2020social}, our model (\ref{eq:2}) is simpler and retains the higher-order coupling information. More significantly, we can discuss the controllability of such systems. In the next subsection, we will establish a Kalman-rank-like condition by exploiting nonlinear control theory. 

\subsection{Controllability of Uniform Hypergraphs}\label{sec:2.4}
If one rewrites the tensor vector multiplications in the multilinear system (\ref{eq:2}) explicitly as in Fig. \ref{fig:3} B, the drift term $\textsf{A}\textbf{x}^{k-1}$ is in fact  a homogeneous polynomial system of degree $k-1$. The controllability of polynomial systems was studied extensively back in 1970s and 80s \cite{doi.org/10.1007/BF01455564,208561,BAILLIEUL1981543,1103422}. In particular, Jurdjevic and Kupka \cite{doi.org/10.1007/BF01455564} obtained strong results in terms of the controllability of homogeneous polynomial systems with constant input multipliers (i.e., $\textbf{b}_j$ are constant vectors). 

\textit{Definition 6 (\cite{doi.org/10.1007/BF01455564}):} A dynamical system is called \textit{strongly controllable} if it can be driven from any initial state to any target state in any instant of positive time given a suitable choice of control inputs. 

\textit{Theorem 1 (\cite{doi.org/10.1007/BF01455564}):} Consider the following system 
\begin{equation}\label{eq:5}
\dot{\textbf{x}} = \textbf{f}(\textbf{x}) + \sum_{j=1}^m\textbf{b}_ju_j.
\end{equation}
Suppose that $\textbf{f}$ is a homogeneous polynomial system of odd degree. Then the system (\ref{eq:5}) is strongly controllable if and only if the rank of the Lie algebra spanned by the set of vector fields $\{\textbf{f}, \textbf{b}_1,\textbf{b}_2,\dots,\textbf{b}_m\}$ is $n$ at all points of $\mathbb{R}^n$.  Moreover, the Lie algebra is of full rank at all points of $\mathbb{R}^n$ if and only if it is of full rank at the origin. 

The rank of the Lie algebra can be found by evaluating the recursive Lie brackets of $\{\textbf{f}, \textbf{b}_1,\textbf{b}_2,\dots,\textbf{b}_m\}$ at the origin. The Lie bracket of two vector fields $\textbf{f}$ and $\textbf{g}$ at a point \textbf{x} is defined as
\begin{equation}
[\textbf{f},\textbf{g}]_\textbf{x} = \nabla\textbf{g}(\textbf{x})\textbf{f}(\textbf{x}) - \nabla\textbf{f}(\textbf{x})\textbf{g}(\textbf{x}),
\end{equation}  
where, $\nabla$ is the gradient operation. Detailed definitions of Lie algebra and Lie brackets can be found in any differential geometry textbook. Based on Theorem 1, we can derive a Kalman-rank-like condition for the tensor-based multilinear system (\ref{eq:2}).

\textit{Definition 7:}
Let $\mathscr{C}_0$ be the linear span of $\{\textbf{b}_1,\textbf{b}_2,\dots,\textbf{b}_m\}$ and $\textsf{A}\in\mathbb{R}^{n\times n\times \dots\times n}$ be a supersymmetric tensor. For each integer $q\geq 1$, define $\mathscr{C}_q$ inductively as the linear span of 
\begin{equation}
\mathscr{C}_{q-1}\cup \{\textsf{A}\textbf{v}_1\textbf{v}_2\dots\textbf{v}_{k-1}|\textbf{v}_l\in \mathscr{C}_{q-1}\}.
\end{equation}
Denote the subspace $\mathscr{C}(\textsf{A}, \textbf{B})=\cup_{q\geq 0} \mathscr{C}_q$ where $\textbf{B}=\begin{bmatrix}\textbf{b}_1 & \textbf{b}_2 & \dots & \textbf{b}_m\end{bmatrix}\in\mathbb{R}^{n\times m}$.

\textit{Corollary 1:}
Suppose that $k$ is even. The multilinear system (\ref{eq:2}) is strongly controllable if and only if the subspace $\mathscr{C}(\textsf{A}, \textbf{B})$ spans $\mathbb{R}^n$, or equivalently, the matrix \textbf{C}, including all the column vectors from $\mathscr{C}(\textsf{A}, \textbf{B})$, has rank $n$. 

\textit{Proof:}
We show that $\mathscr{C}(\textsf{A}, \textbf{B})$ consists of all the recursive Lie brackets of $\{\textsf{A}\textbf{x}^{k-1}, \textbf{b}_1,\textbf{b}_2,\dots,\textbf{b}_m\}$ at the origin. Without loss of generality, assume that $m=1$. Since \textsf{A} is supersymmetric, the recursive Lie brackets are given by (we omit all the scalars in the calculation)
\begin{align*}
[\textbf{b},\textsf{A}\textbf{x}^{k-1}]_{\textbf{0}}&=(\frac{d}{d\textbf{x}}\Bigr|_{\substack{\textbf{x}=\textbf{0}}} \textsf{A}\textbf{x}^{k-1})\textbf{b} = \textbf{0},\\
[\textbf{b}, [\textbf{b}, \textsf{A}\textbf{x}^{k-1}]]_{\textbf{0}} & = (\frac{d}{d\textbf{x}}\Bigr|_{\substack{\textbf{x}=\textbf{0}}}  \textsf{A}\textbf{x}^{k-2}\textbf{b})\textbf{b} = \textbf{0},\\
&\vdots \\
[\textbf{b},[ \dots,[ [\textbf{b}, \textsf{A}\textbf{x}^{k-1}]]]]_{\textbf{0}} & = (\frac{d}{d\textbf{x}}\Bigr|_{\substack{\textbf{x}=\textbf{0}}}  \textsf{A}\textbf{x}\textbf{b}^{k-2})\textbf{b} = \textsf{A}\textbf{b}^{k-1}.
\end{align*}
We then repeat the recursive process for the brackets $[\textsf{A}\textbf{b}^{k-1}, \textsf{A}\textbf{x}^{k-1}]$,  $[\textsf{A}\textbf{b}^{k-1}, \textsf{A}\textbf{x}^{k-2}\textbf{b}]$, $\dots$, $[\textsf{A}\textbf{b}^{k-1}, \textsf{A}\textbf{x}\textbf{b}^{k-2}]$ in the second iteration. After the  $q$-th iteration for some $q$, the subspace $\mathscr{C}(\textsf{A}, \textbf{B})$ contains all the Lie brackets of the vector fields $\{\textsf{A}\textbf{x}^{k-1}, \textbf{b}\}$ at the origin. Lastly, when $k$ is even, the drift term $\textsf{A}\textbf{x}^{k-1}$ is a family of homogeneous polynomial fields of odd degree. Based on Theorem 1, the result follows immediately.
\hfill $\blacksquare$

\textit{Corollary 2:}
Given the subspace $\mathscr{C}(\textsf{A}, \textbf{B})=\cup_{q\geq 0} \mathscr{C}_q$, there exists an integer $q\leq n$ such that $\mathscr{C}(\textsf{A}, \textbf{B}) = \mathscr{C}_q$.

\textit{Proof:}
The proof follows immediately from the fact that $\mathscr{C}(\textsf{A}, \textbf{B})$ is a finite-dimensional vector space \cite{doi.org/10.1007/BF01455564}.
\hfill $\blacksquare$

\begin{algorithm}[ht]
\caption{Computing the reduced controllability matrix.}
\label{alg:1}
\begin{algorithmic}[1]
\STATE{Given a supersymmetric tensor $\textsf{A}\in\mathbb{R}^{n\times n\times \dots \times n}$ and a control matrix $\textbf{B}\in\mathbb{R}^{n\times m}$}\\
\STATE{Unfold \textsf{A} into a matrix \textbf{A} by stacking the last $k-1$ modes, i.e., $\textbf{A}\in\mathbb{R}^{n\times n^{k-1}}$}
\STATE{Set $\textbf{C}_r = \textbf{B}$ and $j=0$}\\
 \WHILE{$j < n$}
\STATE{Compute $\textbf{L}=\textbf{A}(\textbf{C}_r\otimes \textbf{C}_r\otimes \stackrel{k-1}{\cdots} \otimes \textbf{C}_r)$}\\
\STATE{Set $\textbf{C}_r=\begin{bmatrix} \textbf{C}_r & \textbf{L}\end{bmatrix}$}\\
\STATE{Compute the economy-size SVD of $\textbf{C}_r$, and remove the zero singular values, i.e., $\textbf{C}_r=\textbf{U}\textbf{S}\textbf{V}^\top$ where $\textbf{S}\in\mathbb{R}^{s\times s}$, and $s$ is the rank of $\textbf{C}_r$}
\STATE{Set $\textbf{C}_r=\textbf{U}$, and $j=j+1$}
 \ENDWHILE
\RETURN The reduced controllability matrix $\textbf{C}_r$.
\end{algorithmic}
\end{algorithm}

We can treat the matrix \textbf{C} as the controllability matrix of the multilinear system (\ref{eq:2}). When $k=2$ and $q=n-1$, Corollary 1 is reduced to the famous Kalman rank condition for linear systems. However, computing the controllability matrix can be computationally demanding as the number of columns of \textbf{C} grows exponentially with $k$. In Algorithm \ref{alg:1}, we offer a memory-efficient approach to obtain a reduced form of the controllability matrix, denoted by $\textbf{C}_r$, by exploiting the economy-size matrix SVD. Step 2 is referred to as the 1-mode tensor unfolding (matrization) defined in \cite{doi:10.1137/07070111X}. Here $\otimes$ denotes the Kronecker product, and Step 5 can be computed fast without evaluating the actual Kronecker products of $\textbf{C}_r$ \cite{kron}. One may also exploit the sparse tensor/matrix structure to further save the computation and memory.  

\textit{Lemma 1:} Suppose that $\textbf{A}\in\mathbb{R}^{n\times n^{k-1}}$ is defined as in Algorithm 1, and $\textbf{X}\in\mathbb{R}^{n\times m}$ is an arbitrary matrix with rank $s$. Then the following two matrices
\begin{align*}
    \textbf{P}&=\textbf{A}(\textbf{X}\otimes \textbf{X}\otimes \stackrel{k-1}{\cdots} \otimes \textbf{X})\in\mathbb{R}^{n\times m^{k-1}},\\
    \textbf{Q}&=\textbf{A}(\textbf{U}\otimes \textbf{U}\otimes \stackrel{k-1}{\cdots} \otimes \textbf{U})\in\mathbb{R}^{n\times s^{k-1}},
\end{align*}
share the same column space, where $\textbf{U}\in\mathbb{R}^{n\times s}$ is the matrix including the first $s$ left singular vectors of $\textbf{X}$. 

\textit{Proof:} Without loss of generality, assume that $k=3$. Suppose that $\textbf{X}=\textbf{U}\textbf{S}\textbf{V}^\top$ with $\textbf{U}\in\mathbb{R}^{n\times s}$. Based on the properties of the Kronecker product, one can write
\begin{align*}
    \textbf{P}&=\textbf{A}[(\textbf{U}\textbf{S}\textbf{V}^\top)\otimes (\textbf{U}\textbf{S}\textbf{V}^\top)]\\
    & =\textbf{A}[(\textbf{U}\otimes \textbf{U})(\textbf{S}\otimes \textbf{S})(\textbf{V}\otimes \textbf{V})^\top]\\
    & = \textbf{Q}[(\textbf{S}\otimes \textbf{S})(\textbf{V}\otimes \textbf{V})^\top] = \textbf{Q}\tilde{\textbf{S}}\tilde{\textbf{V}}^\top,
\end{align*}
where, $\tilde{\textbf{S}}=\textbf{S}\otimes \textbf{S}\in\mathbb{R}^{s^2\times s^2}$ is a diagonal matrix and $\tilde{\textbf{V}}=\textbf{V}\otimes \textbf{V}\in\mathbb{R}^{m^2\times s^2}$ is a semi-orthogonal matrix. Thus, it follows immediately that $\textbf{P}$ and $\textbf{Q}$ share the same column space. \hfill $\blacksquare$

\textit{Proposition 1:} The column space of the reduced controllability matrix $\textbf{C}_r$ is $\mathscr{C}(\textsf{A}, \textbf{B})$. 

\textit{Proof:} The result follows immediately from Lemma 1 and Corollary 2. \hfill $\blacksquare$

\textit{Remark:} The controllability of homogeneous polynomial systems of even degree is still an open problem to best of authors knowledge \cite{aeyels1984local,melody2003nonlinear}. The reason is intimately related to the fact that the roots of polynomial systems of even degree might all be complex \cite{aeyels1984local}. Therefore, we cannot guarantee the condition for controllability of odd uniform hypergraphs using Corollary 1. Nevertheless, a weaker form of controllability, called \textit{(local) accessibility}, can be obtained for the multilinear system (\ref{eq:2}) with odd $k$ based on the Kalman-rank-like condition. 

Given $\textbf{x}_0\in\mathbb{R}^n$ and control inputs, define $\mathscr{R}(\textbf{x}_0,t)$ to be the set of all $\textbf{x}\in\mathbb{R}^n$ for which the system can be driven from $\textbf{x}_0$ to $\textbf{x}$ at time $t$.

\textit{Definition 8 (\cite{bloch2003nonholonomic}):} A dynamical system is called \textit{accessible} if for any initial state $\textbf{x}_0\in\mathbb{R}^n$ and $T>0$, the reachable set $\mathscr{R}_T(\textbf{x}_0)=\cup_{0\leq t\leq T}\mathscr{R}(\textbf{x}_0,t)$ contains a nonempty open set.

The accessibility of a control system  requires only that the reachable set from a given point contains a nonempty open set, rather than being equal to the whole space $\mathbb{R}^n$ (required for strong controllability). Accessibility holds at a point if the span of the smallest Lie algebra of vector fields containing the drift and input vector fields is $\mathbb{R}^n$ at that point \cite{bloch2003nonholonomic}.

\textit{Corollary 3:}
The multilinear system (\ref{eq:2}) is accessible if the subspace $\mathscr{C}(\textsf{A}, \textbf{B})$ spans $\mathbb{R}^n$, or equivalently, the matrix \textbf{C}, including all the column vectors from $\mathscr{C}(\textsf{A}, \textbf{B})$, has rank $n$. 

\textit{Proof:} The smallest Lie algebra of vector fields containing $\textsf{A}\textbf{x}^{k-1}$ and $\textbf{b}_1,\dots,\textbf{b}_m$ at the origin is $\mathscr{C}(\textsf{A}, \textbf{B})$ by Corollary 1. Based on the second part of Theorem 1 (i.e., the Lie algebra is of full rank at all points of $\mathbb{R}^n$ if and only if it is of full rank at the origin), the result follows immediately. \hfill $\blacksquare$
 
\subsection{MCN of Special Hypergraphs}\label{sec:2.6}

\begin{figure}[t]
\centering
\tcbox[colback=mygray,top=5pt,left=5pt,right=5pt,bottom=5pt]{
\begin{tikzpicture}[scale=0.91, transform shape]
\node[vertex,text=white,scale=0.7] (v1) {1};
\node[vertex,right of=v1,text=white,scale=0.7] (v2) {2};
\node[vertex,right of=v2,text=white,scale=0.7] (v3) {3};
\node[vertex,right of=v3,text=white,scale=0.7] (v4) {4};
\node[vertex,below of=v4,text=white,scale=0.7,mygray] (v5) {};
\node[rectangle,below of=v3,text=black,scale=0.7] (v10) {\Large rank(\textbf{C}) = 4};
\node[vertex,above of=v1,node distance=25pt,mygray,scale=0.7,text=black] (v6) {$b_1$};
\node[vertex,above of=v2,node distance=25pt,mygray,scale=0.7,text=black] (v7) {$b_2$};
\node[vertex,above of=v3,node distance=25pt,mygray,scale=0.7,text=black] (v8) {$b_3$};
\node[rectangle,above right = 0.5pt and 15pt of v4,mygray,scale=0.7,text=black] (v9) {
$\large
\textbf{C} = \begin{bmatrix}
\textbf{b}_1 & \textbf{b}_2 & \textbf{b}_3 & \textsf{A}\textbf{b}_1\textbf{b}_2\textbf{b}_3 & \textsf{A}\textbf{b}_1^3 & \dots 
\end{bmatrix}
$};
\node[rectangle,right of=v5,node distance=90pt,mygray,scale=0.7,text=black] (v9) {
$\Large
=
\begin{bmatrix}
b_1 & 0 & 0 & 0 & 0 & \dots\\
0 & b_2 & 0 & 0 & 0 & \dots \\
0 & 0 & b_3 & 0 & 0 & \dots \\
0 & 0 & 0 & \frac{1}{6}b_1b_2b_3 &0 & \dots
\end{bmatrix}
$};

\path [->,shorten >=1pt,shorten <=1pt, thick](v6) edge node[left] {} (v1);
\path [->,shorten >=1pt,shorten <=1pt, thick](v7) edge node[left] {} (v2);
\path [->,shorten >=1pt,shorten <=1pt, thick](v8) edge node[left] {} (v3);

\begin{pgfonlayer}{background}
\draw[edge,color=orange] (v1) -- (v2) -- (v3) -- (v4);
\end{pgfonlayer}
\end{tikzpicture}
}
\caption{Controllability matrix. 4-uniform hypergraph with four nodes and a hyperedge $\{1,2,3,4\}$, and its  controllability matrix $\textbf{C}$.}
\label{fig:4}
\end{figure}
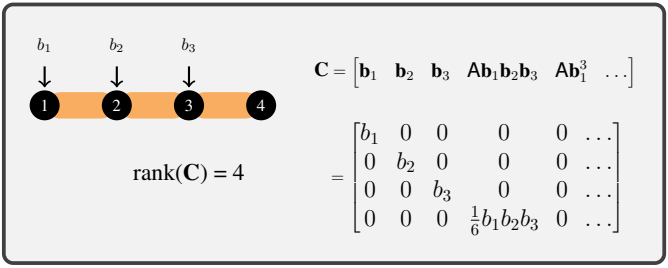

According to Corollary 1 and 2, we can discuss the controllability of even uniform hypergraphs. Similarly to \cite{doi.org/10.1038/nature10011,doi.org/10.1038/ncomms3447}, we want to determine the MCN, denoted by $n^{*}$, whose control is sufficient for achieving controllability of the hypergraph. For example, let's consider the simplest even uniform hypergraph, i.e., the 4-uniform hypergraph with four nodes. We find that the MCN of this hypergraph is three based on the Kalman-rank-like condition, see Fig. \ref{fig:4}. More significantly, we discover that the MCN of even uniform hyperchains, hyperrings and hyperstars as well as complete even uniform hypergraphs behaves similarly to those of chains, rings, stars and complete graphs.  

\textit{Proposition 2:} Suppose that $k$ is even. If \textsf{G} is a $k$-uniform hyperchain with $n$ nodes, then the MCN of \textsf{G} is given by $n^{*}=k-1$.

\textit{Proof:} We first show that when the number of control inputs $m= k-2$, \textsf{G} is never controllable. Based on the definition of tensor vector multiplication and the special structure of the adjacency tensor, it can be shown that for any choice of the control matrix \textbf{B} with $m=k-2$,
\begin{equation*}
    \textsf{A}\textbf{b}_{j_1}\textbf{b}_{j_2}\dots\textbf{b}_{j_{k-1}} = \textbf{0},
\end{equation*}
where, \textsf{A} is the adjacency tensor of \textsf{G}. This is because the set of vectors $\textbf{b}_{j_l}$ must contain a duplicate for one $l=1,2,\dots,k-1$. Therefore, the controllability matrix $\textbf{C}$ has rank $k-2$, and \textsf{G} is not controllable. Next, we present one control strategy with $m= k-1$.

Assume that the first $k-1$ nodes are controlled, and $\textbf{b}_j$ is associated with the node $j$ for $j=1,2,\dots,k-1$. Let
\begin{equation*}
    \textbf{b}_{j}=\textsf{A}\textbf{b}_{j-k+1}\textbf{b}_{j-k+2}\dots\textbf{b}_{j-1}
\end{equation*}
for $j=k,k+1,\dots,n$. Similarly, it can be shown that the following matrix
\begin{equation*}
    \begin{bmatrix} \textbf{b}_1 & \textbf{b}_2 & \dots & \textbf{b}_n\end{bmatrix}\in\mathbb{R}^{n\times n}
\end{equation*} 
is upper triangular with nonzero entries along the diagonal. Hence, the controllability matrix \textbf{C} has rank $n$, and the MCN of \textsf{G} is $k-1$.
\hfill $\blacksquare$

\textit{Proposition 3:} Suppose that $k$ is even and $k\geq 4$. If \textsf{G} is a $k$-uniform hyperring with $n$ nodes and $n>k+1$, then the MCN of \textsf{G} is given by $n^{*}=k-1$.

\textit{Proof:} We first note that when $n=k+1$, \textsf{G} is complete, see Proposition 5 below. For $n>k+1$, the first part of the proof follows in exactly the same fashion as in Proposition 2. Moreover, due to the special structure of the adjacency tensor and the definition of tensor vector multiplication, it is straightforward to show that for $k\geq 4$, controlling the first $k-1$ nodes will be enough to make the rank of the controllability matrix \textbf{C} equal to $n$. Therefore, the MCN of \textsf{G} is $k-1$. \hfill $\blacksquare$

\textit{Proposition 4:} Suppose that $k$ is even. If \textsf{G} is a $k$-uniform hyperstar with $n$ nodes and $n>k$, then the MCN of \textsf{G} is given by $n^{*}=n-2$.

\textit{Proof:} We first note that when $n=k$, \textsf{G} is a uniform hyperchain. For $n>k$, we show that when the number of control inputs $m= n-3$, \textsf{G} is never controllable. Let $m_{\text{int}}$ denote the number of control inputs selected from the set of the internal nodes. Then $m-m_{\text{int}}$ is the number of control inputs selected from the set of the leaf nodes. According to the definition of tensor vector multiplication and the special structure of the adjacency tensor, it can be shown that

\textit{Case 1:} When $m_{\text{int}} = k-1$, $\text{rank}(\textbf{C}) = n-2$;

\textit{Case 2:} When $m_{\text{int}} = k-2$, $\text{rank}(\textbf{C}) = n-1$;

\textit{Case 3:} When $0\leq m_{\text{int}} < k-2$, $\text{rank}(\textbf{C}) = n-3$.\\
Therefore, \textsf{G} is not controllable. However, if one adds one more input from the set of the leaf nodes in Case 2, the rank of the controllability matrix \textbf{C} will reach $n$. Hence, the MCN of \textsf{G} is $n-2$.
\hfill $\blacksquare$

\textit{Proposition 5:} Suppose that $k$ is even. If \textsf{G} is a complete $k$-uniform hypergraph with $n$ nodes, then the MCN of \textsf{G} is given by $n^{*}=n-1$.

\textit{Proof:} We first show that when the number of control inputs $m= n-2$, \textsf{G} is never controllable. Without loss of generality, assume that the first $n-2$ nodes are controlled, and $\textbf{b}_j$ is associated with the node $j$ for $j=1,2,\dots,n-2$. According to the definition of tensor vector multiplication and the fact that \textsf{G} is complete, it can be shown that all the column vectors in the controllability matrix \textbf{C} have the last two entries equal. Thus, the rank of \textbf{C} is equal to $n-1$, and \textsf{G} is not controllable. However, when the first $n-1$ nodes are controlled, the last node can be easily reached by any combination of $(k-1)$ $\textbf{b}_j$ for $j=1,2,\dots,n-1$ since \textsf{G} is complete, and the controllability matrix \textbf{C} has rank $n$. Therefore, the MCN of \textsf{G} is $n-1$.
 \hfill $\blacksquare$

Proposition 2, 4 and 5 are valid when $k=2$ where the MCN of chains, stars and complete graphs are equal to $1$, $n-2$ and $n-1$, respectively. Similarly to rings, the MCN of even uniform hyperrings does not depend on $n$. However, Proposition 3 does not hold for $k = 2$ because the pairwise transitions between nodes would produce linearly dependent relations in the standard rings. Furthermore, we discover that the MCN of odd uniform hyperchains, hyperrings and hyperstars as well as complete odd uniform hypergraphs, identified by the Kalman-rank-like condition, follows exactly the same patterns as stated in Proposition 2 to 5 even though the condition only offers accessibility, see Section \ref{sec:3.1}. Nonetheless, controllability or accessibility of a uniform hypergragh is closely associated to its underlying architecture, and the MCN can be used to  quantify some topological attributes of a uniform hypergraph.

\section{Hypergraph Robustness}\label{sec:2.7}


Network robustness is the ability of a network to survive from random failures or deliberate attacks (e.g., removal of nodes or edges) \cite{barabasi2003linked,callaway2000network,abbas2012robust}. It is intimately related to the underlying network structure/topology. Many measures have been proposed to quantify the robustness of a graph, and one of the popular measures is called effective resistance \cite{spielman2011graph,ELLENS20112491}. We propose MCN as a measure of hypergraph robustness since it can provide insights into the topology and connectivity of uniform hypergraphs according to their controllability or accessibility. Intuitively, if the MCN of a uniform hypergraph is high,  it will take more effort/energy to control the hypergraph or steer the underlying system. 

In Section \ref{sec:3.2} and \ref{sec:3.1}, we will determine the MCN of even and odd uniform hyperchains, hyperrings and hyperstars, and their different variants in simulated datasets. As expected, the rules for selecting the minimum subset of ``control nodes'' of odd uniform hypergraphs associated with the MCN follow the same patterns as these of even uniform hypergraphs. More interestingly, we find that the MCN of these configurations are related to their degree distributions, and the high degree nodes are often preferred as control nodes. 

\subsection{Control Nodes Selection with MCN}

In order to measure the robustness of uniform hypergraphs, we need to efficiently determine their MCN. However, finding the MCN of uniform hypergraphs using a brute-force search will be computationally demanding. We provide a heuristic approach for estimating the minimum subset of control nodes of a uniform hypergraph in which nodes are chosen based on the maximum change in the rank of the reduced controllability matrix, see Algorithm \ref{alg:2}. Here $\textbf{C}_D$ denotes the reduced controllability matrix formed from the inputs in the index set $D$, and can be computed using Algorithm \ref{alg:1}. If a uniform hypergraph is non-connected, we can first identify the connected components (which can be defined similarly as in graphs), and then apply the algorithm to each component. In Step 7, if multiple $s^*$ are obtained, we can pick one randomly, or use some other conditions to break the tie, e.g., by selecting the node with the highest degree. It turns out that Algorithm \ref{alg:2} with high likelihood can find the MCN of a medium-sized uniform hypergraph, and it is much faster than a brute-force search, see Section \ref{sec:3.5}.

\begin{algorithm}[ht]
\caption{Greedy Control Nodes Selection with MCN.}
\label{alg:2}
\begin{algorithmic}[1]
\STATE{Given a supersymmetric tensor $\textsf{A}\in\mathbb{R}^{n\times n\times \dots \times n}$}\\
\STATE{Let $S=\{1,2,\dots,n\}$ and $D=\emptyset$}\\
\WHILE{$\text{rank}(\textbf{C}_D)<n$}
\FOR{$s\in S\setminus D$}
\STATE{Compute $\Delta(s)=\text{rank}(\textbf{C}_{D\cup \{s\}})-\text{rank}(\textbf{C}_D)$ using Algorithm \ref{alg:1}}\\
\ENDFOR
\STATE{Set $s^{*} = \text{argmax}_{s\in S\setminus D}\Delta(s)$}\\
\STATE{Set $D=D\cup \{s^*\}$}
\ENDWHILE
\RETURN A subset of control nodes $D$.
\end{algorithmic}
\end{algorithm}

\section{Numerical Examples}\label{sec:3}
All the numerical examples presented were performed on a Linux machine with 16 GB RAM and a 2.0 GHz Quad-Core Intel Core i5 processor in MATLAB R2020a.

\subsection{Even Uniform Hypergraphs}\label{sec:3.2}
Recall that in $k$-uniform hyperchains, hyperrings and hyperstars, every hyperegde has exactly $k-1$ overlapping nodes. However, for $k\geq 3$, one can relax this requirement allowing number of overlapping nodes between hyperedges to vary, and obtain variants of the three configurations. We consider the case where every intersection between hyperedges contains $r$ nodes for $0<r<k-1$, and denote these variants by $r$-hyperchains, $r$-hyperrings and $r$-hyperstars. Note that $k$-uniform hyperchains, hyperrings and hyperstars are the cases where $r=k-1$. Uniform hyperchains, hyperrings and hyperstars, and their different variants have many applications in real life. For example, we can use them to model complex networks such as computer networks, supply chains and organizational hierarchy. Understanding the control mechanisms of these configurations will be greatly beneficial for achieving network security, efficient communications and energy savings. 

\begin{figure}[t!]
\centering
\tcbox[colback=mygray,top=0pt,left=-2pt,right=0pt,bottom=5pt]{
\begin{tikzpicture}[scale=0.65, transform shape]
\node[vertex,text=white,scale=0.7] (v1) {1};
\node[vertex,right of=v1,text=white,scale=0.7,node distance = 20pt] (v2) {2};
\node[vertex,right of=v2,text=white,scale=0.7,node distance = 20pt] (v3) {3};
\node[vertex,right of=v3,text=white,scale=0.7,node distance = 20pt] (v4) {4};
\node[vertex,right of=v4,text=white,scale=0.7,node distance = 20pt] (v5) {5};
\node[vertex,right of=v5,text=white,scale=0.7,node distance = 20pt] (v6) {6};
\node[vertex,right of=v6,text=white,scale=0.7,node distance = 20pt] (v7) {7};
\node[vertex,left of=v1,node distance=30pt,mygray,text=black,scale=0.7] (text1) {\huge\textbf{A}};

\node[vertex,above of=v2,node distance=20pt,mygray,scale=0.7] (v8) {};
\node[vertex,above of=v3,node distance=20pt,mygray,scale=0.7] (v9) {};
\node[vertex,above of=v4,node distance=20pt,mygray,scale=0.7] (v10) {};
\node[vertex,above of=v5,node distance=20pt,mygray,scale=0.7] (v11) {};
\node[vertex,above of=v6,node distance=20pt,mygray,scale=0.7] (v12) {};

\node[vertex,right of=v7,text=white,scale=0.7,node distance=30pt] (v13) {1};
\node[vertex,right of=v13,text=white,scale=0.7,node distance = 20pt] (v14) {2};
\node[vertex,right of=v14,text=white,scale=0.7,node distance = 20pt] (v15) {3};
\node[vertex,right of=v15,text=white,scale=0.7,node distance = 20pt] (v16) {4};
\node[vertex,right of=v16,text=white,scale=0.7,node distance = 20pt] (v17) {5};
\node[vertex,right of=v17,text=white,scale=0.7,node distance = 20pt] (v18) {6};
\node[vertex,right of=v18,text=white,scale=0.7,node distance = 20pt] (v19) {7};
\node[vertex,right of=v19,text=white,scale=0.7,node distance = 20pt] (v20) {8};
\node[vertex,right of=v20,text=white,scale=0.7,node distance = 20pt] (v21) {9};
\node[vertex,right of=v21,text=white,scale=0.58,node distance = 20pt] (v22) {10};

\node[vertex,above of=v14,node distance=20pt,mygray,scale=0.7] (v23) {};
\node[vertex,above of=v15,node distance=20pt,mygray,scale=0.7] (v24) {};
\node[vertex,above of=v16,node distance=20pt,mygray,scale=0.7] (v25) {};
\node[vertex,above of=v17,node distance=20pt,mygray,scale=0.7] (v26) {};
\node[vertex,above of=v19,node distance=20pt,mygray,scale=0.7] (v27) {};
\node[vertex,above of=v20,node distance=20pt,mygray,scale=0.7] (v28) {};
\node[vertex,above of=v21,node distance=20pt,mygray,scale=0.7] (v29) {};

\path [->,shorten >=1pt,shorten <=1pt, thick](v8) edge node[left] {} (v2);
\path [->,shorten >=1pt,shorten <=1pt, thick](v9) edge node[left] {} (v3);
\path [->,shorten >=1pt,shorten <=1pt, thick,cyan](v10) edge node[left] {} (v4);
\path [->,shorten >=1pt,shorten <=1pt, thick](v11) edge node[left] {} (v5);
\path [->,shorten >=1pt,shorten <=1pt, thick](v12) edge node[left] {} (v6);

\path [->,shorten >=1pt,shorten <=1pt, thick](v23) edge node[left] {} (v14);
\path [->,shorten >=1pt,shorten <=1pt, thick](v24) edge node[left] {} (v15);
\path [->,shorten >=1pt,shorten <=1pt, thick,cyan](v25) edge node[left] {} (v16);
\path [->,shorten >=1pt,shorten <=1pt, thick](v26) edge node[left] {} (v17);
\path [->,shorten >=1pt,shorten <=1pt, thick,cyan](v27) edge node[left] {} (v19);
\path [->,shorten >=1pt,shorten <=1pt, thick](v28) edge node[left] {} (v20);
\path [->,shorten >=1pt,shorten <=1pt, thick](v29) edge node[left] {} (v21);

\node[vertex,below of=v1,text=white,scale=0.7,node distance=30pt] (a1) {1};
\node[vertex,right of=a1,text=white,scale=0.7,node distance = 20pt] (a2) {2};
\node[vertex,right of=a2,text=white,scale=0.7,node distance = 20pt] (a3) {3};
\node[vertex,right of=a3,text=white,scale=0.7,node distance = 20pt] (a4) {4};
\node[vertex,right of=a4,text=white,scale=0.7,node distance = 20pt] (a5) {5};
\node[vertex,right of=a5,text=white,scale=0.7,node distance = 20pt] (a6) {6};
\node[vertex,left of=a1,node distance=30pt,mygray,text=black,scale=0.7] (text2.1) {\huge\textbf{B}};

\node[vertex,above of=a2,node distance=20pt,mygray,scale=0.7] (a15) {};
\node[vertex,above of=a3,node distance=20pt,mygray,scale=0.7] (a16) {};
\node[vertex,above of=a4,node distance=20pt,mygray,scale=0.7] (a17) {};
\node[vertex,above of=a5,node distance=20pt,mygray,scale=0.7] (a18) {};

\node[vertex,below of=v13,text=white,scale=0.7,node distance=30pt] (a7) {1};
\node[vertex,right of=a7,text=white,scale=0.7,node distance = 20pt] (a8) {2};
\node[vertex,right of=a8,text=white,scale=0.7,node distance = 20pt] (a9) {3};
\node[vertex,right of=a9,text=white,scale=0.7,node distance = 20pt] (a10) {4};
\node[vertex,right of=a10,text=white,scale=0.7,node distance = 20pt] (a11) {5};
\node[vertex,right of=a11,text=white,scale=0.7,node distance = 20pt] (a12) {6};
\node[vertex,right of=a12,text=white,scale=0.7,node distance = 20pt] (a13) {7};
\node[vertex,right of=a13,text=white,scale=0.7,node distance = 20pt] (a14) {8};

\node[vertex,above of=a8,node distance=20pt,mygray,scale=0.7] (a19) {};
\node[vertex,above of=a9,node distance=20pt,mygray,scale=0.7] (a20) {};
\node[vertex,above of=a10,node distance=20pt,mygray,scale=0.7] (a21) {};
\node[vertex,above of=a11,node distance=20pt,mygray,scale=0.7] (a22) {};
\node[vertex,above of=a13,node distance=20pt,mygray,scale=0.7] (a23) {};

\path [->,shorten >=1pt,shorten <=1pt, thick](a15) edge node[left] {} (a2);
\path [->,shorten >=1pt,shorten <=1pt, thick,cyan](a16) edge node[left] {} (a3);
\path [->,shorten >=1pt,shorten <=1pt, thick,cyan](a17) edge node[left] {} (a4);
\path [->,shorten >=1pt,shorten <=1pt, thick](a18) edge node[left] {} (a5);

\path [->,shorten >=1pt,shorten <=1pt, thick](a19) edge node[left] {} (a8);
\path [->,shorten >=1pt,shorten <=1pt, thick,cyan](a20) edge node[left] {} (a9);
\path [->,shorten >=1pt,shorten <=1pt, thick,cyan](a21) edge node[left] {} (a10);
\path [->,shorten >=1pt,shorten <=1pt, thick,cyan](a22) edge node[left] {} (a11);
\path [->,shorten >=1pt,shorten <=1pt, thick](a23) edge node[left] {} (a13);

\node[vertex,below of=a1,text=white,scale=0.7,node distance=30pt] (t1) {1};
\node[vertex,right of=t1,text=white,scale=0.7,node distance = 20pt] (t2) {2};
\node[vertex,right of=t2,text=white,scale=0.7,node distance = 20pt] (t3) {3};
\node[vertex,right of=t3,text=white,scale=0.7,node distance = 20pt] (t4) {4};
\node[vertex,right of=t4,text=white,scale=0.7,node distance = 20pt] (t5) {5};
\node[vertex,left of=t1,node distance=30pt,mygray,text=black,scale=0.7] (text2) {\huge\textbf{C}};

\node[vertex,above of=t4,node distance=20pt,mygray,scale=0.7] (t6) {};
\node[vertex,above of=t2,node distance=20pt,mygray,scale=0.7] (t7) {};
\node[vertex,above of=t3,node distance=20pt,mygray,scale=0.7] (t8) {};

\node[vertex,below of=a7,text=white,scale=0.7,node distance=30pt] (t9) {1};
\node[vertex,right of=t9,text=white,scale=0.7,node distance = 20pt] (t10) {2};
\node[vertex,right of=t10,text=white,scale=0.7,node distance = 20pt] (t11) {3};
\node[vertex,right of=t11,text=white,scale=0.7,node distance = 20pt] (t12) {4};
\node[vertex,right of=t12,text=white,scale=0.7,node distance = 20pt] (t13) {5};
\node[vertex,right of=t13,text=white,scale=0.7,node distance = 20pt] (t14) {6};

\node[vertex,above of=t12,node distance=20pt,mygray,scale=0.7] (t15) {};
\node[vertex,above of=t10,node distance=20pt,mygray,scale=0.7] (t16) {};
\node[vertex,above of=t11,node distance=20pt,mygray,scale=0.7] (t17) {};

\path [->,shorten >=1pt,shorten <=1pt, thick,cyan](t6) edge node[left] {} (t4);
\path [->,shorten >=1pt,shorten <=1pt, thick,cyan](t7) edge node[left] {} (t2);
\path [->,shorten >=1pt,shorten <=1pt, thick,cyan](t8) edge node[left] {} (t3);

\path [->,shorten >=1pt,shorten <=1pt, thick,cyan](t15) edge node[left] {} (t12);
\path [->,shorten >=1pt,shorten <=1pt, thick](t16) edge node[left] {} (t10);
\path [->,shorten >=1pt,shorten <=1pt, thick,cyan](t17) edge node[left] {} (t11);

\node[vertex,below of=t1,text=white,scale=0.7,node distance=90pt] (s1) {1};
\node[vertex,above right = 11pt and 1.5pt of s1,text=white,scale=0.7,node distance = 20pt] (s2) {2};
\node[vertex,above right = 11pt and 1.5pt of s2,text=white,scale=0.7,node distance = 20pt] (s3) {3};
\node[vertex,above right = 11pt and 1.5pt of s3,text=white,scale=0.7,node distance = 20pt] (s4) {4};
\node[vertex,below right = 11pt and 1.5pt of s4,text=white,scale=0.7,node distance = 20pt] (s5) {5};
\node[vertex,below right = 11pt and 1.5pt of s5,text=white,scale=0.7,node distance = 20pt] (s6) {6};
\node[vertex,below right = 11pt and 1.5pt of s6,text=white,scale=0.7,node distance = 20pt] (s7) {7};
\node[vertex,left of=s7,text=white,scale=0.7,node distance = 20pt] (s8) {8};
\node[vertex,left of=s8,text=white,scale=0.7,node distance = 20pt] (s9) {9};
\node[vertex,below of=text2,node distance=60pt,mygray,text=black,scale=0.7] (text4) {\huge\textbf{D}};

\node[vertex,above of=s1,node distance=20pt,mygray,scale=0.7] (s10) {};
\node[vertex,above of=s2,node distance=20pt,mygray,scale=0.7] (s11) {};
\node[vertex,above of=s4,node distance=20pt,mygray,scale=0.7] (s12) {};
\node[vertex,above of=s5,node distance=20pt,mygray,scale=0.7] (s13) {};
\node[vertex,above of=s7,node distance=20pt,mygray,scale=0.7] (s14) {};
\node[vertex,above of=s8,node distance=20pt,mygray,scale=0.7] (s15) {};

\node[vertex,below of=t9,text=white,scale=0.7,node distance=90pt] (s16) {1};
\node[vertex,above of=s16,text=white,scale=0.7,node distance = 20pt] (s17) {2};
\node[vertex,above of=s17,text=white,scale=0.7,node distance = 20pt] (s18) {3};
\node[vertex,above of=s18,text=white,scale=0.7,node distance = 20pt] (s19) {4};
\node[vertex,right of=s19,text=white,scale=0.7,node distance = 20pt] (s20) {5};
\node[vertex,right of=s20,text=white,scale=0.7,node distance = 20pt] (s21) {6};
\node[vertex,right of=s21,text=white,scale=0.7,node distance = 20pt] (s22) {7};
\node[vertex,below of=s22,text=white,scale=0.7,node distance = 20pt] (s23) {8};
\node[vertex,below of=s23,text=white,scale=0.7,node distance = 20pt] (s24) {9};
\node[vertex,below of=s24,text=white,scale=0.58,node distance = 20pt] (s25) {10};
\node[vertex,left of=s25,text=white,scale=0.58,node distance = 20pt] (s26) {11};
\node[vertex,left of=s26,text=white,scale=0.58,node distance = 20pt] (s27) {12};

\node[vertex,above of=s19,node distance=20pt,mygray,scale=0.7] (s28) {};
\node[vertex,above of=s22,node distance=20pt,mygray,scale=0.7] (s29) {};
\node[vertex,left of=s16,node distance=20pt,mygray,scale=0.7] (s30) {};
\node[vertex,right of=s25,node distance=20pt,mygray,scale=0.7] (s31) {};
\node[vertex,left of=s17,node distance=20pt,mygray,scale=0.7] (s32) {};
\node[vertex,right of=s23,node distance=20pt,mygray,scale=0.7] (s33) {};
\node[vertex,above of=s20,node distance=20pt,mygray,scale=0.7] (s34) {};
\node[vertex,above of=s26,node distance=20pt,mygray,scale=0.7] (s35) {};

\path [->,shorten >=1pt,shorten <=1pt, thick,cyan](s10) edge node[left] {} (s1);
\path [->,shorten >=1pt,shorten <=1pt, thick](s11) edge node[left] {} (s2);
\path [->,shorten >=1pt,shorten <=1pt, thick,cyan](s12) edge node[left] {} (s4);
\path [->,shorten >=1pt,shorten <=1pt, thick](s13) edge node[left] {} (s5);
\path [->,shorten >=1pt,shorten <=1pt, thick,cyan](s14) edge node[left] {} (s7);
\path [->,shorten >=1pt,shorten <=1pt, thick](s15) edge node[left] {} (s8);

\path [->,shorten >=1pt,shorten <=1pt, thick,cyan](s28) edge node[left] {} (s19);
\path [->,shorten >=1pt,shorten <=1pt, thick,cyan](s29) edge node[left] {} (s22);
\path [->,shorten >=1pt,shorten <=1pt, thick,cyan](s30) edge node[left] {} (s16);
\path [->,shorten >=1pt,shorten <=1pt, thick,cyan](s31) edge node[left] {} (s25);
\path [->,shorten >=1pt,shorten <=1pt, thick](s32) edge node[left] {} (s17);
\path [->,shorten >=1pt,shorten <=1pt, thick](s33) edge node[left] {} (s23);
\path [->,shorten >=1pt,shorten <=1pt, thick](s34) edge node[left] {} (s20);
\path [->,shorten >=1pt,shorten <=1pt, thick](s35) edge node[left] {} (s26);

\node[vertex,below of=s1,text=white,scale=0.7,node distance=50pt] (b1) {1};
\node[vertex,above of=b1,text=white,scale=0.7,node distance = 20pt] (b2) {2};
\node[vertex,right of=b2,text=white,scale=0.7,node distance = 20pt] (b3) {3};
\node[vertex,right of=b3,text=white,scale=0.7,node distance = 20pt] (b4) {4};
\node[vertex,below of=b4,text=white,scale=0.7,node distance = 20pt] (b5) {5};
\node[vertex,left of=b5,text=white,scale=0.7,node distance = 20pt] (b6) {6};
\node[vertex,below of=text4,node distance=70pt,mygray,text=black,scale=0.7] (text5) {\huge\textbf{E}};

\node[vertex,above of=b2,node distance=20pt,mygray,scale=0.7] (b15) {};
\node[vertex,above of=b3,node distance=20pt,mygray,scale=0.7] (b16) {};
\node[vertex,above of=b4,node distance=20pt,mygray,scale=0.7] (b17) {};
\node[vertex,right of=b5,node distance=20pt,mygray,scale=0.7] (b18) {};

\node[vertex,below of=s16,text=white,scale=0.7,node distance=50pt] (b7) {1};
\node[vertex,above of=b7,text=white,scale=0.7,node distance = 20pt] (b8) {2};
\node[vertex,right of=b8,text=white,scale=0.7,node distance = 20pt] (b9) {3};
\node[vertex,right of=b9,text=white,scale=0.7,node distance = 20pt] (b10) {4};
\node[vertex,right of=b10,text=white,scale=0.7,node distance = 20pt] (b11) {5};
\node[vertex,below of=b11,text=white,scale=0.7,node distance = 20pt] (b12) {6};
\node[vertex,left of=b12,text=white,scale=0.7,node distance = 20pt] (b13) {7};
\node[vertex,left of=b13,text=white,scale=0.7,node distance = 20pt] (b14) {8};

\node[vertex,above of=b8,node distance=20pt,mygray,scale=0.7] (b19) {};
\node[vertex,above of=b9,node distance=20pt,mygray,scale=0.7] (b20) {};
\node[vertex,above of=b10,node distance=20pt,mygray,scale=0.7] (b21) {};
\node[vertex,above of=b11,node distance=20pt,mygray,scale=0.7] (b22) {};

\path [->,shorten >=1pt,shorten <=1pt, thick,cyan](b15) edge node[left] {} (b2);
\path [->,shorten >=1pt,shorten <=1pt, thick,cyan](b16) edge node[left] {} (b3);
\path [->,shorten >=1pt,shorten <=1pt, thick,cyan](b17) edge node[left] {} (b4);
\path [->,shorten >=1pt,shorten <=1pt, thick,cyan](b18) edge node[left] {} (b5);

\path [->,shorten >=1pt,shorten <=1pt, thick,cyan](b19) edge node[left] {} (b8);
\path [->,shorten >=1pt,shorten <=1pt, thick,cyan](b20) edge node[left] {} (b9);
\path [->,shorten >=1pt,shorten <=1pt, thick,cyan](b21) edge node[left] {} (b10);
\path [->,shorten >=1pt,shorten <=1pt, thick,cyan](b22) edge node[left] {} (b11);
\path [->,shorten >=1pt,shorten <=1pt, thick,cyan](b10) edge node[left] {} (b13);

\node[vertex,below of=b7,text=white,scale=0.7,node distance=50pt] (q1) {1};
\node[vertex,above of= q1,text=white,scale=0.7,node distance = 20pt] (q2) {2};
\node[vertex,right of=q2,text=white,scale=0.7,node distance = 20pt] (q3) {3};
\node[vertex,right of=q3,text=white,scale=0.7,node distance = 20pt] (q99) {4};
\node[vertex,below right = 1.2pt and 11pt of q99,text=white,scale=0.7,node distance = 20pt] (q4) {5};
\node[vertex,below left = 1.2pt and 11pt of q4,text=white,scale=0.7,node distance = 20pt] (q5) {6};
\node[vertex,left of=q5,text=white,scale=0.7,node distance = 20pt] (q98) {7};
\node[vertex,below of=text5,node distance=50pt,mygray,text=black,scale=0.7] (text6) {\huge\textbf{F}};

\node[vertex,above of=q2,node distance=20pt,mygray,scale=0.7] (q7) {};
\node[vertex,above of=q3,node distance=20pt,mygray,scale=0.7] (q8) {};
\node[vertex,above of=q99,node distance=20pt,mygray,scale=0.7] (q9) {};

\node[vertex,below of=b1,text=white,scale=0.7,node distance=50pt] (q10) {1};
\node[vertex,above of=q10,text=white,scale=0.7,node distance = 20pt] (q11) {2};
\node[vertex,right of=q11,text=white,scale=0.7,node distance = 20pt] (q12) {3};
\node[vertex,right of=q12,text=white,scale=0.7,node distance = 20pt] (q13) {4};
\node[vertex,below of=q13,text=white,scale=0.7,node distance = 20pt] (q14) {5};
\node[vertex,left of=q14,text=white,scale=0.7,node distance = 20pt] (q15) {6};

\node[vertex,above of=q11,node distance=20pt,mygray,scale=0.7] (q17) {};
\node[vertex,above of=q12,node distance=20pt,mygray,scale=0.7] (q18) {};
\node[vertex,above of=q13,node distance=20pt,mygray,scale=0.7] (q19) {};

\path [->,shorten >=1pt,shorten <=1pt, thick,cyan](q7) edge node[left] {} (q2);
\path [->,shorten >=1pt,shorten <=1pt, thick,cyan](q8) edge node[left] {} (q3);
\path [->,shorten >=1pt,shorten <=1pt, thick,cyan](q9) edge node[left] {} (q99);

\path [->,shorten >=1pt,shorten <=1pt, thick,cyan](q17) edge node[left] {} (q11);
\path [->,shorten >=1pt,shorten <=1pt, thick,cyan](q18) edge node[left] {} (q12);
\path [->,shorten >=1pt,shorten <=1pt, thick,cyan](q19) edge node[left] {} (q13);

\node[vertex,below of=q10,text=white,scale=0.7,node distance=90pt] (r1) {1};
\node[vertex,right of=r1,text=white,scale=0.7,node distance = 20pt] (r2) {2};
\node[vertex,right of=r2,text=white,scale=0.7,node distance = 20pt] (r3) {3};
\node[vertex,right of=r3,text=white,scale=0.7,node distance = 20pt] (r4) {4};
\node[vertex,right of= r4,text=white,scale=0.7,node distance = 20pt] (r5) {5};
\node[vertex,right of=r5,text=white,scale=0.7,node distance = 20pt] (r6) {6};
\node[vertex,right of= r6,text=white,scale=0.7,node distance = 20pt] (r7) {7};
\node[vertex,above of= r4,text=white,scale=0.7,node distance = 20pt] (r8) {8};
\node[vertex,above of=r8,text=white,scale=0.7,node distance = 20pt] (r9) {9};
\node[vertex,above of= r9,text=white,scale=0.58,node distance = 20pt] (r10) {10};
\node[vertex,below of=text6,node distance=80pt,mygray,text=black,scale=0.7] (text7) {\huge\textbf{G}};

\node[vertex,above of=r2,node distance=20pt,mygray,scale=0.7] (r11) {};
\node[vertex,above of=r3,node distance=20pt,mygray,scale=0.7] (r12) {};
\node[vertex,above of=r5,node distance=20pt,mygray,scale=0.7] (r13) {};
\node[vertex,above of=r6,node distance=20pt,mygray,scale=0.7] (r14) {};
\node[vertex,right of=r8,node distance=20pt,mygray,scale=0.7] (r15) {};
\node[vertex,right of=r9,node distance=20pt,mygray,scale=0.7] (r16) {};
\node[vertex,above left = 7pt and 7pt of r4,node distance=20pt,mygray,scale=0.7] (r17) {};

\node[vertex,below of=q1,text=white,scale=0.7,node distance=90pt] (r18) {1};
\node[vertex,right of=r18,text=white,scale=0.7,node distance = 20pt] (r19) {2};
\node[vertex,right of=r19,text=white,scale=0.7,node distance = 20pt] (r20) {3};
\node[vertex,right of=r20,text=white,scale=0.7,node distance = 20pt] (r21) {4};
\node[vertex,right of=r21,text=white,scale=0.7,node distance = 20pt] (r22) {5};
\node[vertex,right of=r22,text=white,scale=0.7,node distance = 20pt] (r23) {6};
\node[vertex,right of=r23,text=white,scale=0.7,node distance = 20pt] (r24) {7};
\node[vertex,above right = 11pt and 0.5pt of r21,text=white,scale=0.7,node distance = 20pt] (r25) {8};
\node[vertex,above right = 11pt and 0.5pt of r25,text=white,scale=0.7,node distance = 20pt] (r26) {9};
\node[vertex,above right = 11pt and 0.5pt of r26,text=white,scale=0.58,node distance = 20pt] (r27) {10};
\node[vertex,above left = 11pt and 0.5pt of r21,text=white,scale=0.58,node distance = 20pt] (r28) {11};
\node[vertex,above left = 11pt and 0.5pt of r28,text=white,scale=0.58,node distance = 20pt] (r29) {12};
\node[vertex,above left = 11pt and 0.5pt of r29,text=white,scale=0.58,node distance = 20pt] (r30) {13};

\node[vertex,above of=r19,node distance=20pt,mygray,scale=0.7] (r31) {};
\node[vertex,above of=r20,node distance=20pt,mygray,scale=0.7] (r32) {};
\node[vertex,above of=r22,node distance=20pt,mygray,scale=0.7] (r33) {};
\node[vertex,above of=r23,node distance=20pt,mygray,scale=0.7] (r34) {};
\node[vertex,right of=r25,node distance=20pt,mygray,scale=0.7] (r35) {};
\node[vertex,right of=r26,node distance=20pt,mygray,scale=0.7] (r36) {};
\node[vertex,left of=r28,node distance=20pt,mygray,scale=0.7] (r37) {};
\node[vertex,left of=r29,node distance=20pt,mygray,scale=0.7] (r38) {};
\node[vertex,above of=r21,node distance=20pt,mygray,scale=0.7] (r39) {};

\path [->,shorten >=1pt,shorten <=1pt, thick](r11) edge node[left] {} (r2);
\path [->,shorten >=1pt,shorten <=1pt, thick](r12) edge node[left] {} (r3);
\path [->,shorten >=1pt,shorten <=1pt, thick](r13) edge node[left] {} (r5);
\path [->,shorten >=1pt,shorten <=1pt, thick](r14) edge node[left] {} (r6);
\path [->,shorten >=1pt,shorten <=1pt, thick](r15) edge node[left] {} (r8);
\path [->,shorten >=1pt,shorten <=1pt, thick](r16) edge node[left] {} (r9);
\path [->,shorten >=1pt,shorten <=1pt, thick,cyan](r17) edge node[left] {} (r4);

\path [->,shorten >=1pt,shorten <=1pt, thick](r31) edge node[left] {} (r19);
\path [->,shorten >=1pt,shorten <=1pt, thick](r32) edge node[left] {} (r20);
\path [->,shorten >=1pt,shorten <=1pt, thick](r33) edge node[left] {} (r22);
\path [->,shorten >=1pt,shorten <=1pt, thick](r34) edge node[left] {} (r23);
\path [->,shorten >=1pt,shorten <=1pt, thick](r35) edge node[left] {} (r25);
\path [->,shorten >=1pt,shorten <=1pt, thick](r36) edge node[left] {} (r26);
\path [->,shorten >=1pt,shorten <=1pt, thick](r37) edge node[left] {} (r28);
\path [->,shorten >=1pt,shorten <=1pt, thick](r38) edge node[left] {} (r29);
\path [->,shorten >=1pt,shorten <=1pt, thick,cyan](r39) edge node[left] {} (r21);

\node[vertex,below of=r1,text=white,scale=0.7,node distance=70pt] (c1) {1};
\node[vertex,right of=c1,text=white,scale=0.7,node distance = 20pt] (c2) {2};
\node[vertex,right of=c2,text=white,scale=0.7,node distance = 20pt] (c3) {3};
\node[vertex,right of=c3,text=white,scale=0.7,node distance = 20pt] (c4) {4};
\node[vertex,right of= c4,text=white,scale=0.7,node distance = 20pt] (c5) {5};
\node[vertex,right of=c5,text=white,scale=0.7,node distance = 20pt] (c6) {6};
\node[vertex,above of= c4,text=white,scale=0.7,node distance = 20pt] (c7) {7};
\node[vertex,above of=c7,text=white,scale=0.7,node distance = 20pt] (c8) {8};
\node[vertex,below of=text7,node distance=70pt,mygray,text=black,scale=0.7] (text8) {\huge\textbf{H}};

\node[vertex,above of=c2,node distance=20pt,mygray,scale=0.7] (c19) {};
\node[vertex,above of=c3,node distance=20pt,mygray,scale=0.7] (c20) {};
\node[vertex,above of=c5,node distance=20pt,mygray,scale=0.7] (c21) {};
\node[vertex,right of=c7,node distance=20pt,mygray,scale=0.7] (c22) {};
\node[vertex,above left = 7pt and 7pt of c4,node distance=20pt,mygray,scale=0.7] (c23) {};

\node[vertex,below of=r18,text=white,scale=0.7,node distance=70pt] (c9) {1};
\node[vertex,right of=c9,text=white,scale=0.7,node distance = 20pt] (c10) {2};
\node[vertex,right of=c10,text=white,scale=0.7,node distance = 20pt] (c11) {3};
\node[vertex,right of=c11,text=white,scale=0.7,node distance = 20pt] (c12) {4};
\node[vertex,right of= c12,text=white,scale=0.7,node distance = 20pt] (c13) {5};
\node[vertex,right of=c13,text=white,scale=0.7,node distance = 20pt] (c14) {6};
\node[vertex,above of= c11,text=white,scale=0.7,node distance = 20pt] (c15) {7};
\node[vertex,above of=c15,text=white,scale=0.7,node distance = 20pt] (c16) {8};
\node[vertex,above of= c12,text=white,scale=0.7,node distance = 20pt] (c17) {9};
\node[vertex,above of=c17,text=white,scale=0.58,node distance = 20pt] (c18) {10};

\node[vertex,above of=c10,node distance=20pt,mygray,scale=0.7] (c24) {};
\node[vertex,above of=c13,node distance=20pt,mygray,scale=0.7] (c25) {};
\node[vertex,left of=c15,node distance=20pt,mygray,scale=0.7] (c26) {};
\node[vertex,right of=c17,node distance=20pt,mygray,scale=0.7] (c27) {};
\node[vertex,above left = 7pt and 7pt of c11,node distance=20pt,mygray,scale=0.7] (c28) {};
\node[vertex,above right = 7pt and 7pt of c12,node distance=20pt,mygray,scale=0.7] (c29) {};

\path [->,shorten >=1pt,shorten <=1pt, thick](c19) edge node[left] {} (c2);
\path [->,shorten >=1pt,shorten <=1pt, thick,cyan](c20) edge node[left] {} (c3);
\path [->,shorten >=1pt,shorten <=1pt, thick](c21) edge node[left] {} (c5);
\path [->,shorten >=1pt,shorten <=1pt, thick](c22) edge node[left] {} (c7);
\path [->,shorten >=1pt,shorten <=1pt, thick,cyan](c23) edge node[left] {} (c4);

\path [->,shorten >=1pt,shorten <=1pt, thick](c24) edge node[left] {} (c10);
\path [->,shorten >=1pt,shorten <=1pt, thick](c25) edge node[left] {} (c13);
\path [->,shorten >=1pt,shorten <=1pt, thick](c26) edge node[left] {} (c15);
\path [->,shorten >=1pt,shorten <=1pt, thick](c27) edge node[left] {} (c17);
\path [->,shorten >=1pt,shorten <=1pt, thick,cyan](c28) edge node[left] {} (c11);
\path [->,shorten >=1pt,shorten <=1pt, thick,cyan](c29) edge node[left] {} (c12);

\node[vertex,below of=c1,text=white,scale=0.7,node distance=50pt] (f1) {1};
\node[vertex,right of=f1,text=white,scale=0.7,node distance=20pt] (f2) {2};
\node[vertex,right of=f2,text=white,scale=0.7,node distance=20pt] (f3) {3};
\node[vertex,right of=f3,text=white,scale=0.7,node distance=20pt] (f4) {4};
\node[vertex,below of= f4,text=white,scale=0.7,node distance=20pt] (f5) {5};
\node[vertex,above of=f4,text=white,scale=0.7,node distance=20pt] (f6) {6};
\node[vertex,below of=text8,node distance=70pt,mygray,text=black,scale=0.7] (text6) {\huge\textbf{I}};

\node[vertex,above of=f1,node distance=20pt,mygray,scale=0.7] (f7) {};
\node[vertex,above of=f2,node distance=20pt,mygray,scale=0.7] (f8) {};
\node[vertex,above of=f3,node distance=20pt,mygray,scale=0.7] (f9) {};
\node[vertex,above of=f6,node distance=20pt,mygray,scale=0.7] (f10) {};

\node[vertex,below of=c9,text=white,scale=0.7,node distance=50pt] (f11) {1};
\node[vertex,right of=f11,text=white,scale=0.7,node distance=20pt] (f12) {2};
\node[vertex,right of=f12,text=white,scale=0.7,node distance=20pt] (f13) {3};
\node[vertex,right of=f13,text=white,scale=0.7,node distance=20pt] (f14) {4};
\node[vertex,right of= f14,text=white,scale=0.7,node distance=20pt] (f15) {5};
\node[vertex,above of=f14,text=white,scale=0.7,node distance=20pt] (f16) {6};
\node[vertex,below of=f14,text=white,scale=0.7,node distance=20pt] (f17) {7};

\node[vertex,above of=f11,node distance=20pt,mygray,scale=0.7] (f18) {};
\node[vertex,above of=f12,node distance=20pt,mygray,scale=0.7] (f19) {};
\node[vertex,above of=f13,node distance=20pt,mygray,scale=0.7] (f20) {};
\node[vertex,above of=f15,node distance=20pt,mygray,scale=0.7] (f21) {};
\node[vertex,above of=f16,node distance=20pt,mygray,scale=0.7] (f22) {};

\path [->,shorten >=1pt,shorten <=1pt, thick](f7) edge node[left] {} (f1);
\path [->,shorten >=1pt,shorten <=1pt, thick,cyan](f8) edge node[left] {} (f2);
\path [->,shorten >=1pt,shorten <=1pt, thick,cyan](f9) edge node[left] {} (f3);
\path [->,shorten >=1pt,shorten <=1pt, thick](f10) edge node[left] {} (f6);

\path [->,shorten >=1pt,shorten <=1pt, thick](f18) edge node[left] {} (f11);
\path [->,shorten >=1pt,shorten <=1pt, thick,cyan](f19) edge node[left] {} (f12);
\path [->,shorten >=1pt,shorten <=1pt, thick,cyan](f20) edge node[left] {} (f13);
\path [->,shorten >=1pt,shorten <=1pt, thick](f21) edge node[left] {} (f15);
\path [->,shorten >=1pt,shorten <=1pt, thick](f22) edge node[left] {} (f16);

\begin{pgfonlayer}{background}
\begin{scope}[transparency group,opacity=.9]
\draw[edge,color=orange, line width=8pt] (v1) -- (v2) -- (v3) -- (v4);
\draw[edge,color=red, line width=8pt] (v4) -- (v5) -- (v6) -- (v7);

\draw[edge,color=orange,line width=8pt] (v13) -- (v14) -- (v15) -- (v16);
\draw[edge,color=red,line width=8pt] (v16) -- (v17) -- (v18) -- (v19);
\draw[edge,color=green,line width=8pt] (v19) -- (v20) -- (v21) -- (v22);

\draw[edge,color=orange, line width=8pt] (a1) -- (a2) -- (a3) -- (a4);
\draw[edge,color=red,line width=6pt] (a3) -- (a4) -- (a5) -- (a6);

\draw[edge,color=orange,line width=8pt] (a7) -- (a8) -- (a9) -- (a10);
\draw[edge,color=red,line width=6pt] (a9) -- (a10) -- (a11) -- (a12);
\draw[edge,color=green,line width=8pt] (a11) -- (a12) -- (a13) -- (a14);

\draw[edge,color=orange,line width=8pt] (t1) -- (t2) -- (t3) -- (t4);
\draw[edge,color=red,line width=6pt] (t2) -- (t3) -- (t4) -- (t5);

\draw[edge,color=orange,line width=8pt] (t9) -- (t10) -- (t11) -- (t12);
\draw[edge,color=red,line width=6pt] (t10) -- (t11) -- (t12) -- (t13);
\draw[edge,color=green,line width=8pt] (t11) -- (t12) -- (t13) -- (t14);

\draw[edge,color=orange,line width=8pt] (s1) -- (s2) -- (s3) -- (s4);
\draw[edge,color=red,line width=8pt] (s4) -- (s5) -- (s6) -- (s7);
\draw[edge,color=green,line width=8pt] (s7) -- (s8) -- (s9) -- (s1);

\draw[edge,color=orange,line width=8pt] (s16) -- (s17) -- (s18) -- (s19);
\draw[edge,color=red,line width=8pt] (s19) -- (s20) -- (s21) -- (s22);
\draw[edge,color=green,line width=8pt] (s22) -- (s23) -- (s24) -- (s25);
\draw[edge,color=blue,line width=8pt] (s25) -- (s26) -- (s27) -- (s16);

\draw[edge,color=orange,line width=8pt] (b1) -- (b2) -- (b3) -- (b4);
\draw[edge,color=red,line width=6pt] (b3) -- (b4) -- (b5) -- (b6);
\draw[edge,color=green,line width=10pt] (b5) -- (b6) -- (b1) -- (b2);

\draw[edge,color=orange,line width=8pt] (b7) -- (b8) -- (b9) -- (b10);
\draw[edge,color=red,line width=6pt] (b9) -- (b10) -- (b11) -- (b12);
\draw[edge,color=green,line width=8pt] (b11) -- (b12) -- (b13) -- (b14);
\draw[edge,color=blue,line width=10pt] (b13) -- (b14) -- (b7) -- (b8);

\draw[edge,color=orange,line width=8pt] (q1) -- (q2) -- (q3) -- (q99);
\draw[edge,color=red,line width=6pt] (q2) -- (q3) -- (q99) -- (q4);
\draw[edge,color=green,line width=10pt] (q3) -- (q99) -- (q4) -- (q5);
\draw[edge,color=blue,line width=12pt] (q99) -- (q4) -- (q5) -- (q98);
\draw[edge,color=yellow,line width=4pt] (q4) -- (q5) -- (q98) -- (q1);
\draw[edge,color=pink,line width=2pt] (q5) -- (q98) -- (q1) -- (q2);
\draw[edge,color=brown,line width=10pt] (q98) -- (q1) -- (q2) -- (q3);

\draw[edge,color=orange,line width=8pt] (q10) -- (q11) -- (q12) -- (q13);
\draw[edge,color=red,line width=6pt] (q11) -- (q12) -- (q13) -- (q14);
\draw[edge,color=green,line width=10pt] (q12) -- (q13) -- (q14) -- (q15);
\draw[edge,color=blue,line width=12pt] (q13) -- (q14) -- (q15) -- (q10);
\draw[edge,color=yellow,line width=4pt] (q14) -- (q15) -- (q10) -- (q11);
\draw[edge,color=pink,line width=2pt] (q15) -- (q10) -- (q11) -- (q12);

\draw[edge,color=orange,line width=8pt] (r1) -- (r2) -- (r3) -- (r4);
\draw[edge,color=red,line width=8pt] (r4) -- (r5) -- (r6) -- (r7);
\draw[edge,color=green,line width=8pt] (r4) -- (r8) -- (r9) -- (r10);

\draw[edge,color=orange,line width=8pt] (r18) -- (r19) -- (r20) -- (r21);
\draw[edge,color=red,line width=8pt] (r21) -- (r22) -- (r23) -- (r24);
\draw[edge,color=green,line width=8pt] (r21) -- (r25) -- (r26) -- (r27);
\draw[edge,color=blue,line width=8pt] (r21) -- (r28) -- (r29) -- (r30);

\draw[edge,color=orange,line width=8pt] (c1) -- (c2) -- (c3) -- (c4);
\draw[edge,color=red,line width=6pt] (c3) -- (c4) -- (c5) -- (c6);
\draw[edge,color=green,line width=10pt] (c3) -- (c4) -- (c7) -- (c8);

\draw[edge,color=orange,line width=8pt] (c9) -- (c10) -- (c11) -- (c12);
\draw[edge,color=red,line width=6pt] (c11) -- (c12) -- (c13) -- (c14);
\draw[edge,color=green,line width=10pt] (c12) -- (c11) -- (c15) -- (c16);
\draw[edge,color=blue,line width=4pt] (c11) -- (c12) -- (c17) -- (c18);

\draw[edge,color=orange,line width=8pt] (f1) -- (f2) -- (f3) -- (f4);
\draw[edge,color=red,line width=6pt] (f2) -- (f3) -- (f4) -- (f5);
\draw[edge,color=green,line width=10pt] (f2) -- (f3) -- (f4) -- (f6);

\draw[edge,color=orange,line width=8pt] (f11) -- (f12) -- (f13) -- (f14);
\draw[edge,color=red,line width=12pt] (f12) -- (f13) -- (f14) -- (f15);
\draw[edge,color=green,line width=10pt] (f12) -- (f13) -- (f14) -- (f16);
\draw[edge,color=blue,line width=6pt] (f12) -- (f13) -- (f14) -- (f17);
\end{scope}
\end{pgfonlayer}
\end{tikzpicture}
}
\caption{MCN of 4-uniform hyperchains, hyperrings and hyperstars, and their variants. (A), (B) and (C): 4-uniform 1-hyperchains, 2-hyperchains and hyperchains. (D), (E) and (F): 4-uniform 1-hyperrings, 2-hyperrings and hyperrings. (G), (H) and (I): 4-uniform 1-hyperstars, 2-hyperstars and hyperstars. The nodes with arrows are denoted as the control nodes, and the cyan arrows indicate the control nodes with highest degrees in the configurations.}
\label{fig:22}
\end{figure}
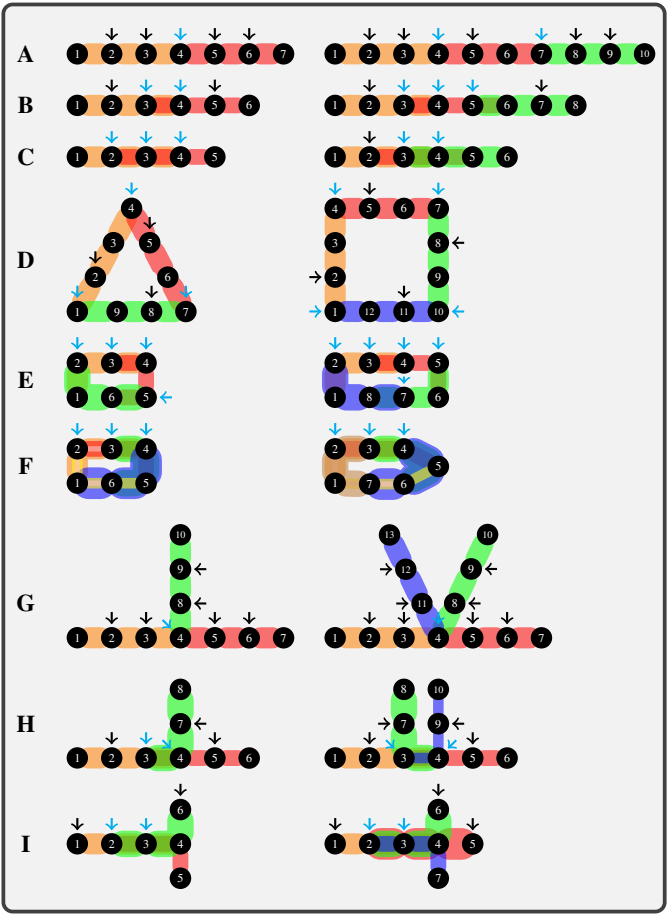

In this example, we determine the MCN of 4-uniform hyperchains, hyperrings and hyperstars, and their variants. The results are shown in Fig. \ref{fig:22}. We only present the most representative minimum subset of control nodes for each configuration. First, the MCN of 4-uniform hyperchains, hyperrings and hyperstars is consistent with the results stated in Proposition 2 to 4. Controlling 4-uniform hyperchains and hyperrings only requires control of three nodes, and controlling 4-uniform hyperstars requires control of $n-2$ nodes, see Fig. \ref{fig:22} C, F and I. Moreover, we discover that the MCN of these configurations is related to their degree distributions. Intuitively, controlling the high degree nodes is the easiest and most natural way for achieving hypergraph control. In particular, all the hypergraph configurations in Fig. \ref{fig:22} contain at least one control node with the highest degree in the corresponding degree distributions. For 4-uniform 1-hyperchains, 1-hyperrings and 1-hyperstars, in which there is one common node between hyperedges, the MCN can be achieved when all the 2-degree nodes are controlled with each hyperedge having three control nodes, see Fig. \ref{fig:22} A, D and G. Furthermore, the control strategies for 4-uniform 2-hyperchains, 2-hyperrings and 2-hyperstars, in which there are two common nodes between hyperedges, are more like some combinations of the previous two, which also require controlling nodes with the highest degree, see Fig. \ref{fig:22} B, E and H. However, it is possible that low degree nodes can accomplish the same goal. For example, the minimum subset of control nodes $\{1,2,3,5,6,8,9\}$  can also achieve the full control of the 1-hyperchain with ten nodes. 

We summarize the MCN of 4-uniform 1-, 2-hyperchains, 1-,2-hyperrings and 1-,2-hyperstars with $n$ nodes in Table \ref{tab:1} (column $n^{*}_{\text{4-unif}}$) based on our observations (i.e., they are not proved). Interestingly, for cases of 2-hyperchains, 2-hyperings and 2-hyperstars, we find that the MCN is $\frac{n+2}{2}$. Whether this results hold in general needs to be further investigated. Moreover, one may easily obtain the MCN of some hybrids of uniform hyperchians, hyperrings and hyperstars, and their variants according to the control strategies discussed above.

\begin{table}[t]
\caption{MCN of the variants of 4- and 3-uniform hyperchains, hyperrings and hyperstars based on our observations. Note that 3-uniform 2-hyperchains, 2-hyperrings and 2-hyperstars are the 3-uniform hyperchains, hyperrings and hyperstars.}
\centering
\begin{tabular}{|l|l|l|}
\hline
\textbf{Configuration}  & $\boldsymbol{n}_{\text{4-unif}}^{*}$ & $\boldsymbol{n}_{\text{3-unif}}^{*}$\\ \hline
1-Hyperchain                                & $\frac{2n+1}{3}$ & $\frac{n+1}{2}$  \\ \hline
1-Hyperring                                  & $\frac{2n}{3}$  & $\frac{n}{2}$    \\ \hline
1-Hyperstar                                 & $\frac{2n+1}{3}$ & $\frac{n+1}{2}$  \\ \hline
2-Hyperchain                                      & $\frac{n+2}{2}$   & $2$         \\ \hline
2-Hyperring                                            & $\frac{n+2}{2}$   & $2$   \\ \hline
2-Hyperstar                                             & $\frac{n+2}{2}$   & $n-2$   \\ \hline
\end{tabular}
\label{tab:1}
\end{table}

\subsection{Odd Uniform Hypergraphs}\label{sec:3.1}

\begin{figure}[t]
\centering
\tcbox[colback=mygray,top=0pt,left=0pt,right=5pt,bottom=5pt]{
\begin{tikzpicture}[scale=0.92, transform shape]
\node[vertex,text=white,scale=0.7] (v1) {1};
\node[vertex,right of=v1,node distance = 20pt,text=white,scale=0.7] (v2) {2};
\node[vertex,right of=v2,node distance = 20pt,text=white,scale=0.7] (v3) {3};
\node[vertex,right of=v3,node distance = 20pt,text=white,scale=0.7] (v4) {4};
\node[vertex,right of=v4,node distance = 20pt,text=white,scale=0.7] (v5) {5};

\node[vertex,above of=v2,node distance=20pt,mygray,scale=0.7] (w1) {};
\node[vertex,above of=v3,node distance=20pt,mygray,scale=0.7] (w2) {};
\node[vertex,above of=v4,node distance=20pt,mygray,scale=0.7] (w3) {};
\path [->,shorten >=1pt,shorten <=1pt, thick](w1) edge node[left] {} (v2);
\path [->,shorten >=1pt,shorten <=1pt, thick,cyan](w2) edge node[left] {} (v3);
\path [->,shorten >=1pt,shorten <=1pt, thick](w3) edge node[left] {} (v4);

\node[vertex, right of=v5, node distance = 20pt,text=white,scale=0.7] (v6) {1};
\node[vertex,right of=v6,node distance = 20pt,text=white,scale=0.7] (v7) {2};
\node[vertex,right of=v7,node distance = 20pt,text=white,scale=0.7] (v8) {3};
\node[vertex,right of=v8,node distance = 20pt,text=white,scale=0.7] (v9) {4};
\node[vertex,right of=v9,node distance = 20pt,text=white,scale=0.7] (v10) {5};
\node[vertex,right of=v10,node distance = 20pt,text=white,scale=0.7] (v11) {6};
\node[vertex,right of=v11,node distance = 20pt,text=white,scale=0.7] (v12) {7};

\node[vertex,above of=v7,node distance=20pt,mygray,scale=0.7] (w4) {};
\node[vertex,above of=v8,node distance=20pt,mygray,scale=0.7] (w5) {};
\node[vertex,above of=v10,node distance=20pt,mygray,scale=0.7] (w6) {};
\node[vertex,above of=v11,node distance=20pt,mygray,scale=0.7] (w7) {};
\path [->,shorten >=1pt,shorten <=1pt, thick](w4) edge node[left] {} (v7);
\path [->,shorten >=1pt,shorten <=1pt, thick,cyan](w5) edge node[left] {} (v8);
\path [->,shorten >=1pt,shorten <=1pt, thick,cyan](w6) edge node[left] {} (v10);
\path [->,shorten >=1pt,shorten <=1pt, thick](w7) edge node[left] {} (v11);

\node[vertex,left of=v1, node distance = 20pt,mygray, text=black,scale=1] (a) {\textbf{A}};

\node[vertex,below of=v1,node distance = 30pt,text=white,scale=0.7] (v13) {1};
\node[vertex,right of=v13,node distance = 20pt,text=white,scale=0.7] (v14) {2};
\node[vertex,right of=v14,node distance = 20pt,text=white,scale=0.7] (v15) {3};
\node[vertex,right of=v15,node distance = 20pt,text=white,scale=0.7] (v16) {4};

\node[vertex,above of=v15,node distance=20pt,mygray,scale=0.7] (w8) {};
\node[vertex,above of=v14,node distance=20pt,mygray,scale=0.7] (w9) {};
\path [->,shorten >=1pt,shorten <=1pt, thick,cyan](w8) edge node[left] {} (v15);
\path [->,shorten >=1pt,shorten <=1pt, thick,cyan](w9) edge node[left] {} (v14);

\node[vertex,below of=v6,node distance = 30pt,text=white,scale=0.7] (v17) {1};
\node[vertex,right of=v17,node distance = 20pt,text=white,scale=0.7] (v18) {2};
\node[vertex,right of=v18,node distance = 20pt,text=white,scale=0.7] (v19) {3};
\node[vertex,right of=v19,node distance = 20pt,text=white,scale=0.7] (v20) {4};
\node[vertex,right of=v20,node distance = 20pt,text=white,scale=0.7] (v21) {5};

\node[vertex,above of=v19,node distance=20pt,mygray,scale=0.7] (w10) {};
\node[vertex,above of=v18,node distance=20pt,mygray,scale=0.7] (w11) {};
\path [->,shorten >=1pt,shorten <=1pt, thick, cyan](w10) edge node[left] {} (v19);
\path [->,shorten >=1pt,shorten <=1pt, thick](w11) edge node[left] {} (v18);

\node[vertex,left of=v13, node distance = 20pt,mygray, text=black,scale=1] (b) {\textbf{B}};

\node[vertex,below of=v13,text=white,scale=0.7,node distance=70pt] (v22) {1};
\node[vertex,above right = 11pt and 1.5pt of v22,text=white,scale=0.7] (v23) {2};
\node[vertex,above right = 11pt and 1.5pt of v23,text=white,scale=0.7] (v24) {3};
\node[vertex,below right = 11pt and 1.5pt of v24,text=white,scale=0.7] (v25) {4};
\node[vertex,below right = 11pt and 1.5pt of v25,text=white,scale=0.7] (v26) {5};
\node[vertex,left of=v26,text=white,scale=0.7, node distance = 20pt] (v27) {6};

\node[vertex,above of=v22,node distance=20pt,mygray,scale=0.7] (w12) {};
\node[vertex,above of=v24,node distance=20pt,mygray,scale=0.7] (w13) {};
\node[vertex,above of=v26,node distance=20pt,mygray,scale=0.7] (w14) {};
\path [->,shorten >=1pt,shorten <=1pt, thick,cyan](w12) edge node[left] {} (v22);
\path [->,shorten >=1pt,shorten <=1pt, thick,cyan](w13) edge node[left] {} (v24);
\path [->,shorten >=1pt,shorten <=1pt, thick,cyan](w14) edge node[left] {} (v26);

\node[vertex,below of=v17,text=white,scale=0.7,node distance=70pt] (v28) {1};
\node[vertex,above of=v28,text=white,scale=0.7,node distance=20pt] (v29) {2};
\node[vertex,above of=v29,text=white,scale=0.7,node distance=20pt] (v30) {3};
\node[vertex,right of=v30,text=white,scale=0.7,node distance=20pt] (v31) {4};
\node[vertex,right of=v31,text=white,scale=0.7,node distance=20pt] (v32) {5};
\node[vertex,below of=v32,text=white,scale=0.7,node distance=20pt] (v33) {6};
\node[vertex,below of=v33,text=white,scale=0.7,node distance=20pt] (v34) {7};
\node[vertex,left of=v34,text=white,scale=0.7,node distance=20pt] (v35) {8};

\node[vertex,left of=v28,node distance=20pt,mygray,scale=0.7] (w15) {};
\node[vertex,above of=v30,node distance=20pt,mygray,scale=0.7] (w16) {};
\node[vertex,above of=v32,node distance=20pt,mygray,scale=0.7] (w17) {};
\node[vertex,right of=v34,node distance=20pt,mygray,scale=0.7] (w18) {};
\path [->,shorten >=1pt,shorten <=1pt, thick,cyan](w15) edge node[left] {} (v28);
\path [->,shorten >=1pt,shorten <=1pt, thick,cyan](w16) edge node[left] {} (v30);
\path [->,shorten >=1pt,shorten <=1pt, thick,cyan](w17) edge node[left] {} (v32);
\path [->,shorten >=1pt,shorten <=1pt, thick,cyan](w18) edge node[left] {} (v34);

\node[vertex,left of=v23, node distance = 30pt,mygray, text=black,scale=1] (c) {\textbf{C}};

\node[vertex,below of=v22,text=white,scale=0.7,node distance=50pt] (v36) {1};
\node[vertex,above of= v36,text=white,scale=0.7,node distance=20pt] (v37) {2};
\node[vertex,right of=v37,text=white,scale=0.7,node distance=20pt] (v38) {3};
\node[vertex,below right = 1.2pt and 13pt of v38,text=white,scale=0.7] (v39) {4};
\node[vertex,below left = 1.2pt and 13pt of v39,text=white,scale=0.7] (v40) {5};

\node[vertex,above of=v37,node distance=20pt,mygray,scale=0.7] (w19) {};
\node[vertex,above of=v38,node distance=20pt,mygray,scale=0.7] (w20) {};
\path [->,shorten >=1pt,shorten <=1pt, thick,cyan](w19) edge node[left] {} (v37);
\path [->,shorten >=1pt,shorten <=1pt, thick,cyan](w20) edge node[left] {} (v38);

\node[vertex,below of=v28,text=white,scale=0.7,node distance=50pt] (v41) {1};
\node[vertex,above of=v41,text=white,scale=0.7,node distance=20pt] (v42) {2};
\node[vertex,right of=v42,text=white,scale=0.7,node distance=20pt] (v43) {3};
\node[vertex,right of=v43,text=white,scale=0.7,node distance=20pt] (v44) {4};
\node[vertex,below of=v44,text=white,scale=0.7,node distance=20pt] (v45) {5};
\node[vertex,left of=v45,text=white,scale=0.7,node distance=20pt] (v46) {6};

\node[vertex,above of=v42,node distance=20pt,mygray,scale=0.7] (w21) {};
\node[vertex,above of=v43,node distance=20pt,mygray,scale=0.7] (w22) {};
\path [->,shorten >=1pt,shorten <=1pt, thick,cyan](w21) edge node[left] {} (v42);
\path [->,shorten >=1pt,shorten <=1pt, thick,cyan](w22) edge node[left] {} (v43);

\node[vertex,left of=v39, node distance = 62pt,mygray, text=black,scale=1] (d) {\textbf{D}};

\node[vertex,below of=v36,text=white,scale=0.7,node distance=70pt] (v47) {1};
\node[vertex,right of=v47,text=white,scale=0.7,node distance = 20pt] (v48) {2};
\node[vertex,right of=v48,text=white,scale=0.7,node distance = 20pt] (v49) {3};
\node[vertex,above right = 11pt and 0.5pt of v49,text=white,scale=0.7] (v50) {4};
\node[vertex,above right = 11pt and 0.5pt of v50,text=white,scale=0.7] (v51) {5};
\node[vertex,below right = 11pt and 0.5pt of v49,text=white,scale=0.7] (v52) {6};
\node[vertex,below right = 11pt and 0.5pt of v52,text=white,scale=0.7] (v53) {7};

\node[vertex,above of=v48,node distance=20pt,mygray,scale=0.7] (w23) {};
\node[vertex,above of=v49,node distance=20pt,mygray,scale=0.7] (w24) {};
\node[vertex,above of=v50,node distance=20pt,mygray,scale=0.7] (w25) {};
\node[vertex,above of=v52,node distance=20pt,mygray,scale=0.7] (w26) {};
\path [->,shorten >=1pt,shorten <=1pt, thick](w23) edge node[left] {} (v48);
\path [->,shorten >=1pt,shorten <=1pt, thick,cyan](w24) edge node[left] {} (v49);
\path [->,shorten >=1pt,shorten <=1pt, thick](w25) edge node[left] {} (v50);
\path [->,shorten >=1pt,shorten <=1pt, thick](w26) edge node[left] {} (v52);

\node[vertex,below of=v41,text=white,scale=0.7,node distance=70pt] (v54) {1};
\node[vertex,right of=v54,text=white,scale=0.7,node distance=20pt] (v55) {2};
\node[vertex,right of=v55,text=white,scale=0.7,node distance=20pt] (v56) {3};
\node[vertex,right of=v56,text=white,scale=0.7,node distance=20pt] (v57) {4};
\node[vertex,right of=v57,text=white,scale=0.7,node distance=20pt] (v58) {5};
\node[vertex,above of=v56,text=white,scale=0.7,node distance=20pt] (v59) {6};
\node[vertex,above of=v59,text=white,scale=0.7,node distance=20pt] (v60) {7};
\node[vertex,below of=v56,text=white,scale=0.7,node distance=20pt] (v61) {8};
\node[vertex,below of=v61,text=white,scale=0.7,node distance=20pt] (v62) {9};

\node[vertex,above of=v55,node distance=20pt,mygray,scale=0.7] (w27) {};
\node[vertex,left of=v61,node distance=20pt,mygray,scale=0.7] (w28) {};
\node[vertex,below of=v57,node distance=20pt,mygray,scale=0.7] (w29) {};
\node[vertex,right of=v59,node distance=20pt,mygray,scale=0.7] (w30) {};
\node[vertex,above left = 7pt and 7pt of v56,node distance=20pt,mygray,scale=0.7] (w31) {};
\path [->,shorten >=1pt,shorten <=1pt, thick](w27) edge node[left] {} (v55);
\path [->,shorten >=1pt,shorten <=1pt, thick](w28) edge node[left] {} (v61);
\path [->,shorten >=1pt,shorten <=1pt, thick](w29) edge node[left] {} (v57);
\path [->,shorten >=1pt,shorten <=1pt, thick](w30) edge node[left] {} (v59);
\path [->,shorten >=1pt,shorten <=1pt, thick,cyan](w31) edge node[left] {} (v56);

\node[vertex,left of=v47, node distance = 20pt,mygray, text=black,scale=1] (e) {\textbf{E}};

\node[vertex,below of=v47,text=white,scale=0.7,node distance=90pt] (v63) {1};
\node[vertex,right of=v63,text=white,scale=0.7,node distance=20pt] (v64) {2};
\node[vertex,right of=v64,text=white,scale=0.7,node distance=20pt] (v65) {3};
\node[vertex,above of=v65,text=white,scale=0.7,node distance=20pt] (v66) {4};
\node[vertex,below of=v65,text=white,scale=0.7,node distance=20pt] (v67) {5};

\node[vertex,above of=v63,node distance=20pt,mygray,scale=0.7] (w32) {};
\node[vertex,above of=v64,node distance=20pt,mygray,scale=0.7] (w33) {};
\node[vertex,above of=v66,node distance=20pt,mygray,scale=0.7] (w34) {};
\path [->,shorten >=1pt,shorten <=1pt, thick](w32) edge node[left] {} (v63);
\path [->,shorten >=1pt,shorten <=1pt, thick,cyan](w33) edge node[left] {} (v64);
\path [->,shorten >=1pt,shorten <=1pt, thick](w34) edge node[left] {} (v66);

\node[vertex,below of=v54,text=white,scale=0.7,node distance=90pt] (v68) {1};
\node[vertex,right of=v68,text=white,scale=0.7,node distance=20pt] (v69) {2};
\node[vertex,right of=v69,text=white,scale=0.7,node distance=20pt] (v70) {3};
\node[vertex,right of=v70,text=white,scale=0.7,node distance=20pt] (v71) {4};
\node[vertex,above of=v70,text=white,scale=0.7,node distance=20pt] (v72) {5};
\node[vertex,below of=v70,text=white,scale=0.7,node distance=20pt] (v73) {6};

\node[vertex,above of=v68,node distance=20pt,mygray,scale=0.7] (w35) {};
\node[vertex,above of=v69,node distance=20pt,mygray,scale=0.7] (w36) {};
\node[vertex,above of=v71,node distance=20pt,mygray,scale=0.7] (w37) {};
\node[vertex,above of=v72,node distance=20pt,mygray,scale=0.7] (w38) {};
\path [->,shorten >=1pt,shorten <=1pt, thick](w35) edge node[left] {} (v68);
\path [->,shorten >=1pt,shorten <=1pt, thick,cyan](w36) edge node[left] {} (v69);
\path [->,shorten >=1pt,shorten <=1pt, thick](w37) edge node[left] {} (v71);
\path [->,shorten >=1pt,shorten <=1pt, thick](w38) edge node[left] {} (v72);

\node[vertex,left of=v63, node distance = 20pt,mygray, text=black,scale=1] (f) {\textbf{F}};

\begin{pgfonlayer}{background}
\draw[edge,color=orange] (v1) -- (v2) -- (v3);
\draw[edge,color=red] (v3) -- (v4) -- (v5);

\draw[edge,color=orange] (v6) -- (v7) -- (v8);
\draw[edge,color=red] (v8) -- (v9) -- (v10);
\draw[edge,color=green] (v10) -- (v11) -- (v12);

\draw[edge,color=orange] (v13) -- (v14) -- (v15);
\draw[edge,color=red,line width=8] (v14) -- (v15) -- (v16);

\draw[edge,color=orange] (v17) -- (v18) -- (v19);
\draw[edge,color=red, line width=8] (v18) -- (v19) -- (v20);
\draw[edge,color=green] (v19) -- (v20) -- (v21);

\draw[edge,color=orange] (v22) -- (v23) -- (v24);
\draw[edge,color=red] (v24) -- (v25) -- (v26);
\draw[edge,color=green] (v26) -- (v27) -- (v22);

\draw[edge,color=orange] (v28) -- (v29) -- (v30);
\draw[edge,color=red] (v30) -- (v31) -- (v32);
\draw[edge,color=green] (v32) -- (v33) -- (v34);
\draw[edge,color=blue] (v34) -- (v35) -- (v28);

\draw[edge,color=orange] (v36) -- (v37) -- (v38);
\draw[edge,color=red,line width=8] (v37) -- (v38) -- (v39);
\draw[edge,color=green] (v38) -- (v39) -- (v40);
\draw[edge,color=blue,line width=8] (v39) -- (v40) -- (v36);
\draw[edge,color=yellow,line width=12] (v40) -- (v36) -- (v37);

\draw[edge,color=orange] (v41) -- (v42) -- (v43);
\draw[edge,color=red,line width=8] (v42) -- (v43) -- (v44);
\draw[edge,color=green] (v43) -- (v44) -- (v45);
\draw[edge,color=blue,line width=8] (v44) -- (v45) -- (v46);
\draw[edge,color=yellow] (v45) -- (v46) -- (v41);
\draw[edge,color=pink,line width=8] (v46) -- (v41) -- (v42);

\draw[edge,color=orange] (v22) -- (v23) -- (v24);
\draw[edge,color=red] (v24) -- (v25) -- (v26);
\draw[edge,color=green] (v26) -- (v27) -- (v22);

\draw[edge,color=orange] (v47) -- (v48) -- (v49);
\draw[edge,color=red] (v49) -- (v50) -- (v51);
\draw[edge,color=green] (v49) -- (v52) -- (v53);

\draw[edge,color=orange] (v54) -- (v55) -- (v56);
\draw[edge,color=red] (v56) -- (v57) -- (v58);
\draw[edge,color=green] (v56) -- (v59) -- (v60);
\draw[edge,color=blue] (v56) -- (v61) -- (v62);

\draw[edge,color=orange] (v63) -- (v64) -- (v65);
\draw[edge,color=red, line width=8] (v64) -- (v65) -- (v66);
\draw[edge,color=green,line width=12] (v64) -- (v65) -- (v67);

\draw[edge,color=orange] (v68) -- (v69) -- (v70);
\draw[edge,color=red, line width=8] (v69) -- (v70) -- (v71);
\draw[edge,color=green,line width=12] (v69) -- (v70) -- (v72);
\draw[edge,color=blue,line width=6] (v69) -- (v70) -- (v73);

\end{pgfonlayer}
\end{tikzpicture}
}
\caption{MCN of 3-uniform hyperchains, hyperrings and hyperstars, and their variants. (A) and (B): 3-uniform 1-hyperchains and hyperchains. (C) and (D): 3-uniform 1-hyperrings and hyperrings. (E) and (F): 3-uniform 1-hyperstars and hyperstars. The nodes with arrows are denoted as the control nodes, and the cyan arrows indicate the control nodes with the highest degrees in the configurations.}
\label{fig:11}
\end{figure}
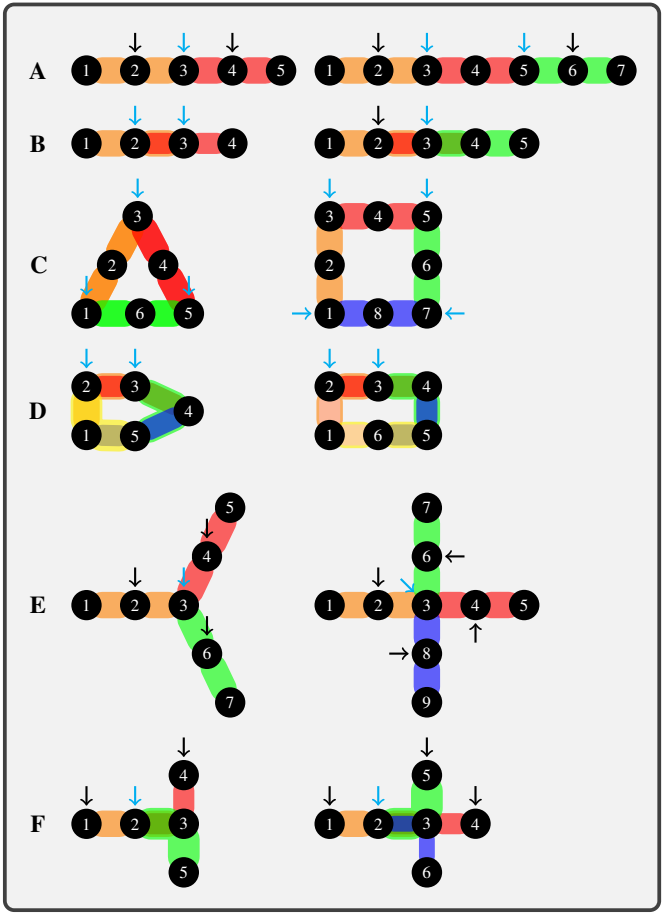

The goal of this example is to determine the MCN of 3-uniform hyperchains, hyperrings and hyperstars, and their variants using the Kalman-rank-like condition even though the condition only offers accessibility for these configurations. The results are shown in Fig. \ref{fig:11}. The rules for selecting the minimum subset of ``control nodes'' of 3-uniform hyperchains, hyperrings and hyperstars, and their variants follow the same patterns as these discussed in Section \ref{sec:3.2}. The MCN of 3-uniform hyperchains, hyperrings and hyperstars are matched with the results stated in Proposition 2 to 4 despite $k$ being odd, see Fig. \ref{fig:11} B, D and F.  Moreover, for 3-uniform 1-hyperchains, 1-hyperrings and 1-hyperstars, in which there is one common node between hyperedges, the MCN can be achieved when all the 2-degree nodes are ``controlled'' with each hyperedge having two ``control nodes'', see Fig. \ref{fig:11} A, C and E. Similarly, we can conclude that the MCN of these configurations are related to their degree distributions. We summarize the MCN of 3-uniform hyperchains, hyperrings and hyperstars, and their variants with $n$ nodes in Table \ref{tab:1} (column $n^{*}_{\text{3-unif}}$) based on our observations. Again, we want to remark that although the ``control nodes'' of 3-uniform hypergraphs may not have physical interpretations in terms of controllability, the MCN can be used to measure hypergraph robustness and detect structural changes, as shown in the following example.

\subsection{Mouse Neuron Endomicroscopy}\label{sec:3.3}
The mouse endomicroscopy dataset is an imaging video created under 10-minute periods of feeding, fasting and re-feeding using fluorescence across space and time in a mouse hypothalamus \cite{9119161,patrick2020}. Twenty neurons are recorded with individual levels of ``firing''. Similar to \cite{9119161}, we want to quantitatively differentiate the three phases using 3-uniform hypergraphs with MCN. First, we compute the multi-correlation among every three neurons, which is defined by
\begin{equation}\label{eq:99}
    \rho = (1-\det{(\textbf{R}}))^{\frac{1}{2}},
\end{equation}
where, $\textbf{R}\in\mathbb{R}^{3\times 3}$ is the correlation matrix of three neuron activity levels \cite{wang2014measures}. When the multi-correlation $\rho$ is greater than a prescribed threshold, we build hyperedges among the three neurons.

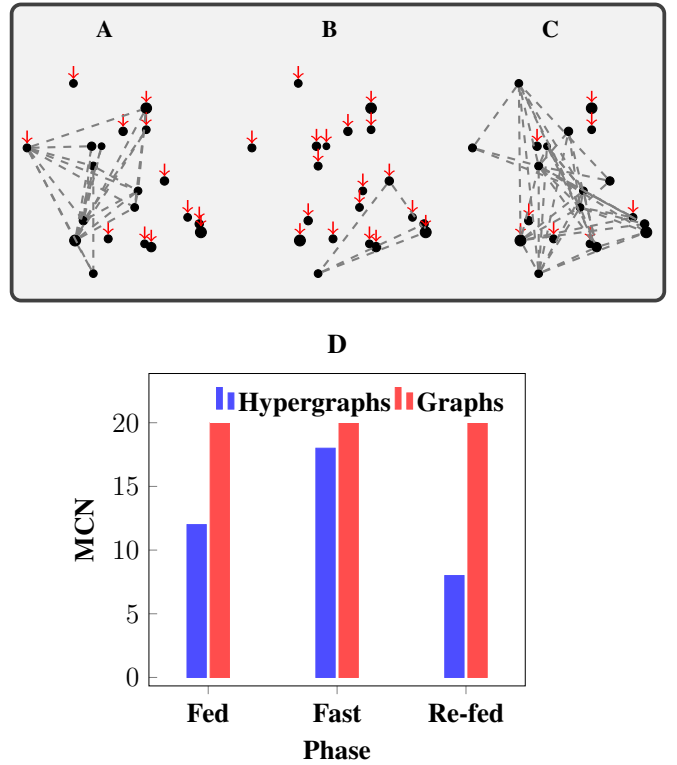
\begin{figure}[t]
    \centering
    \tcbox[colback=mygray,top=0pt,left=-1pt,right=1pt,bottom=5pt]{
    \begin{tikzpicture}[scale=0.022]
    \node[circle,inner sep=0.8] at  (60,150) {\small\textbf{A}};
    
    \node[circle,draw,fill=black,color=black,inner sep=4.12311/4] (v1) at  (110,36) {};
    \node[circle,draw,fill=black,color=black,inner sep=4.24264/4] (v2) at  (13,78) {};
    \node[circle,draw,fill=black,color=black,inner sep=4.47214/4] (v3) at  (96,58) {};
    \node[circle,draw,fill=black,color=black,inner sep=5.65685/4] (v4) at  (85,102) {};
    \node[circle,draw,fill=black,color=black,inner sep=5.83095/4] (v5) at  (42,22) {};
    \node[circle,draw,fill=black,color=black,inner sep=4.47214/4] (v6) at  (47,34) {};
    \node[circle,draw,fill=black,color=black,inner sep=4.24264/4] (v7) at  (53,67) {};
    \node[circle,draw,fill=black,color=black,inner sep=4.12311/4] (v8) at  (85,89) {};
    \node[circle,draw,fill=black,color=black,inner sep=4.12311/4] (v9) at  (78,42) {};
    \node[circle,draw,fill=black,color=black,inner sep=4.12311/4] (v10) at  (53,2) {};
    \node[circle,draw,fill=black,color=black,inner sep=4.24264/4] (v11) at  (62,23) {};
    \node[circle,draw,fill=black,color=black,inner sep=3.60555/4] (v12) at  (58,79) {};
    \node[circle,draw,fill=black,color=black,inner sep=4.47214/4] (v13) at  (52,79) {};
    \node[circle,draw,fill=black,color=black,inner sep=4.47214/4] (v14) at  (71,88) {};
    \node[circle,draw,fill=black,color=black,inner sep=4.12311/4] (v15) at  (80,52) {};
    \node[circle,draw,fill=black,color=black,inner sep=4.47214/4] (v16) at  (117,32) {};
    \node[circle,draw,fill=black,color=black,inner sep=4.12311/4] (v17) at  (41,117) {};
    \node[circle,draw,fill=black,color=black,inner sep=5/4] (v18) at  (88,18) {};
    \node[circle,draw,fill=black,color=black,inner sep=5.83095/4] (v19) at  (118,27) {};
    \node[circle,draw,fill=black,color=black,inner sep=4.12311/4] (v20) at  (84,20) {};
    
    \node[vertex,above of=v1,node distance=10pt,scale=0.7, white,opacity=0] (w1) {};
    \node[vertex,above of=v3,node distance=10pt,scale=0.7, white,opacity=0] (w3) {};
    \node[vertex,above of=v11,node distance=10pt,scale=0.7, white,opacity=0] (w11) {};
    \node[vertex,above of=v14,node distance=10pt,scale=0.7, white,opacity=0] (w14) {};
    \node[vertex,above of=v16,node distance=10pt,scale=0.7, white,opacity=0] (w16) {};
    \node[vertex,above of=v17,node distance=10pt,scale=0.7, white,opacity=0] (w17) {};
    \node[vertex,above of=v18,node distance=10pt,scale=0.7, white,opacity=0] (w18) {};
    \node[vertex,above of=v19,node distance=10pt,scale=0.7, white,opacity=0] (w19) {};
    \node[vertex,above of=v20,node distance=10pt,scale=0.7, white,opacity=0] (w20) {};
    \node[vertex,above of=v2,node distance=10pt,scale=0.7, white,opacity=0] (w2) {};
    \node[vertex,above of=v4,node distance=10pt,scale=0.7, white,opacity=0] (w4) {};
    \node[vertex,above of=v8,node distance=10pt,scale=0.7, white,opacity=0] (w8) {};
    
    \draw [->, red, line width=0.6pt] (w1) -- (v1);
    \draw [->, red, line width=0.6pt] (w3) -- (v3);
    \draw [->, red, line width=0.6pt] (w11) -- (v11);
    \draw [->, red, line width=0.6pt] (w14) -- (v14);
    \draw [->, red, line width=0.6pt] (w16) -- (v16);
    \draw [->, red, line width=0.6pt] (w17) -- (v17);
    \draw [->, red, line width=0.6pt] (w18) -- (v18);
    \draw [->, red, line width=0.6pt] (w19) -- (v19);
    \draw [->, red, line width=0.6pt] (w20) -- (v20);
    \draw [->, red, line width=0.6pt] (w2) -- (v2);
    \draw [->, red, line width=0.6pt] (w4) -- (v4);
    \draw [->, red, line width=0.6pt] (w8) -- (v8);
     
    \draw [dashed, thick,gray] (v2) -- (v4);
    \draw [dashed, thick,gray] (v2) -- (v5);
    \draw [dashed, thick,gray] (v2) -- (v6);
    \draw [dashed, thick,gray] (v2) -- (v7);
    \draw [dashed, thick,gray] (v2) -- (v9);
    \draw [dashed, thick,gray] (v2) -- (v10);
    \draw [dashed, thick,gray] (v2) -- (v13);
    \draw [dashed, thick,gray] (v2) -- (v15);
    \draw [dashed, thick,gray] (v4) -- (v5);
    \draw [dashed, thick,gray] (v4) -- (v9);
    \draw [dashed, thick,gray] (v4) -- (v15);
    \draw [dashed, thick,gray] (v5) -- (v6);
    \draw [dashed, thick,gray] (v5) -- (v7);
    \draw [dashed, thick,gray] (v5) -- (v8);
    \draw [dashed, thick,gray] (v5) -- (v9);
    \draw [dashed, thick,gray] (v5) -- (v10);
    \draw [dashed, thick,gray] (v5) -- (v12);
    \draw [dashed, thick,gray] (v5) -- (v13);
    \draw [dashed, thick,gray] (v5) -- (v15);
    \draw [dashed, thick,gray] (v6) -- (v7);
    \draw [dashed, thick,gray] (v6) -- (v8);
    \draw [dashed, thick,gray] (v6) -- (v12);
    \draw [dashed, thick,gray] (v6) -- (v13);
    \draw [dashed, thick,gray] (v6) -- (v15);
    \draw [dashed, thick,gray] (v7) -- (v10);
    \draw [dashed, thick,gray] (v7) -- (v15);
    \draw [dashed, thick,gray] (v9) -- (v15);

     \end{tikzpicture}
     \hspace{0.2cm}
     \begin{tikzpicture}[scale=0.022]
     \node[circle,inner sep=0.8] at  (60,150) {\small\textbf{B}};
    \node[circle,draw,fill=black,color=black,inner sep=4.12311/4] (v1) at  (110,36) {};
    \node[circle,draw,fill=black,color=black,inner sep=4.24264/4] (v2) at  (13,78) {};
    \node[circle,draw,fill=black,color=black,inner sep=4.47214/4] (v3) at  (96,58) {};
    \node[circle,draw,fill=black,color=black,inner sep=5.65685/4] (v4) at  (85,102) {};
    \node[circle,draw,fill=black,color=black,inner sep=5.83095/4] (v5) at  (42,22) {};
    \node[circle,draw,fill=black,color=black,inner sep=4.47214/4] (v6) at  (47,34) {};
    \node[circle,draw,fill=black,color=black,inner sep=4.24264/4] (v7) at  (53,67) {};
    \node[circle,draw,fill=black,color=black,inner sep=4.12311/4] (v8) at  (85,89) {};
    \node[circle,draw,fill=black,color=black,inner sep=4.12311/4] (v9) at  (78,42) {};
    \node[circle,draw,fill=black,color=black,inner sep=4.12311/4] (v10) at  (53,2) {};
    \node[circle,draw,fill=black,color=black,inner sep=4.24264/4] (v11) at  (62,23) {};
    \node[circle,draw,fill=black,color=black,inner sep=3.60555/4] (v12) at  (58,79) {};
    \node[circle,draw,fill=black,color=black,inner sep=4.47214/4] (v13) at  (52,79) {};
    \node[circle,draw,fill=black,color=black,inner sep=4.47214/4] (v14) at  (71,88) {};
    \node[circle,draw,fill=black,color=black,inner sep=4.12311/4] (v15) at  (80,52) {};
    \node[circle,draw,fill=black,color=black,inner sep=4.47214/4] (v16) at  (117,32) {};
    \node[circle,draw,fill=black,color=black,inner sep=4.12311/4] (v17) at  (41,117) {};
    \node[circle,draw,fill=black,color=black,inner sep=5/4] (v18) at  (88,18) {};
    \node[circle,draw,fill=black,color=black,inner sep=5.83095/4] (v19) at  (118,27) {};
    \node[circle,draw,fill=black,color=black,inner sep=4.12311/4] (v20) at  (84,20) {};
    
    \node[vertex,above of=v1,node distance=10pt,scale=0.7, white,opacity=0] (w1) {};
    \node[vertex,above of=v2,node distance=10pt,scale=0.7, white,opacity=0] (w2) {};
    \node[vertex,above of=v4,node distance=10pt,scale=0.7, white,opacity=0] (w4) {};
    \node[vertex,above of=v5,node distance=10pt,scale=0.7, white,opacity=0] (w5) {};
    \node[vertex,above of=v6,node distance=10pt,scale=0.7, white,opacity=0] (w6) {};
    \node[vertex,above of=v7,node distance=10pt,scale=0.7, white,opacity=0] (w7) {};
    \node[vertex,above of=v8,node distance=10pt,scale=0.7, white,opacity=0] (w8) {};
    \node[vertex,above of=v9,node distance=10pt,scale=0.7, white,opacity=0] (w9) {};
    \node[vertex,above of=v11,node distance=10pt,scale=0.7, white,opacity=0] (w11) {};
    \node[vertex,above of=v12,node distance=10pt,scale=0.7, white,opacity=0] (w12) {};
    \node[vertex,above of=v13,node distance=10pt,scale=0.7, white,opacity=0] (w13) {};
    \node[vertex,above of=v14,node distance=10pt,scale=0.7, white,opacity=0] (w14) {};
    \node[vertex,above of=v15,node distance=10pt,scale=0.7, white,opacity=0] (w15) {};
    \node[vertex,above of=v17,node distance=10pt,scale=0.7, white,opacity=0] (w17) {};
    \node[vertex,above of=v18,node distance=10pt,scale=0.7, white,opacity=0] (w18) {};
    \node[vertex,above of=v20,node distance=10pt,scale=0.7, white,opacity=0] (w20) {};
    
    \node[vertex,above of=v19,node distance=10pt,scale=0.7, white,opacity=0] (w19) {};
    \node[vertex,above of=v3,node distance=10pt,scale=0.7, white,opacity=0] (w3) {};
    
    \draw [->, red, line width=0.6pt] (w1) -- (v1);
    \draw [->, red, line width=0.6pt] (w2) -- (v2);
    \draw [->, red, line width=0.6pt] (w4) -- (v4);
    \draw [->, red, line width=0.6pt] (w5) -- (v5);
    \draw [->, red, line width=0.6pt] (w6) -- (v6);
    \draw [->, red, line width=0.6pt] (w7) -- (v7);
    \draw [->, red, line width=0.6pt] (w8) -- (v8);
    \draw [->, red, line width=0.6pt] (w9) -- (v9);
    \draw [->, red, line width=0.6pt] (w11) -- (v11);
    \draw [->, red, line width=0.6pt] (w12) -- (v12);
    \draw [->, red, line width=0.6pt] (w13) -- (v13);
    \draw [->, red, line width=0.6pt] (w14) -- (v14);
    \draw [->, red, line width=0.6pt] (w15) -- (v15);
    \draw [->, red, line width=0.6pt] (w17) -- (v17);
    \draw [->, red, line width=0.6pt] (w18) -- (v18);
    \draw [->, red, line width=0.6pt] (w20) -- (v20);
    \draw [->, red, line width=0.6pt] (w19) -- (v19);
    \draw [->, red, line width=0.6pt] (w3) -- (v3);

    \draw [dashed, thick,gray] (v3) -- (v10);
    \draw [dashed, thick,gray] (v3) -- (v19);
    \draw [dashed, thick,gray] (v10) -- (v16);
    \draw [dashed, thick,gray] (v10) -- (v19);
    \draw [dashed, thick,gray] (v16) -- (v19);

     \end{tikzpicture}
     \hspace{0.2cm}
     \begin{tikzpicture}[scale=0.022]
     \node[circle,inner sep=0.8] at  (60,150) {\small\textbf{C}};
    \node[circle,draw,fill=black,color=black,inner sep=4.12311/4] (v1) at  (110,36) {};
    \node[circle,draw,fill=black,color=black,inner sep=4.24264/4] (v2) at  (13,78) {};
    \node[circle,draw,fill=black,color=black,inner sep=4.47214/4] (v3) at  (96,58) {};
    \node[circle,draw,fill=black,color=black,inner sep=5.65685/4] (v4) at  (85,102) {};
    \node[circle,draw,fill=black,color=black,inner sep=5.83095/4] (v5) at  (42,22) {};
    \node[circle,draw,fill=black,color=black,inner sep=4.47214/4] (v6) at  (47,34) {};
    \node[circle,draw,fill=black,color=black,inner sep=4.24264/4] (v7) at  (53,67) {};
    \node[circle,draw,fill=black,color=black,inner sep=4.12311/4] (v8) at  (85,89) {};
    \node[circle,draw,fill=black,color=black,inner sep=4.12311/4] (v9) at  (78,42) {};
    \node[circle,draw,fill=black,color=black,inner sep=4.12311/4] (v10) at  (53,2) {};
    \node[circle,draw,fill=black,color=black,inner sep=4.24264/4] (v11) at  (62,23) {};
    \node[circle,draw,fill=black,color=black,inner sep=3.60555/4] (v12) at  (58,79) {};
    \node[circle,draw,fill=black,color=black,inner sep=4.47214/4] (v13) at  (52,79) {};
    \node[circle,draw,fill=black,color=black,inner sep=4.47214/4] (v14) at  (71,88) {};
    \node[circle,draw,fill=black,color=black,inner sep=4.12311/4] (v15) at  (80,52) {};
    \node[circle,draw,fill=black,color=black,inner sep=4.47214/4] (v16) at  (117,32) {};
    \node[circle,draw,fill=black,color=black,inner sep=4.12311/4] (v17) at  (41,117) {};
    \node[circle,draw,fill=black,color=black,inner sep=5/4] (v18) at  (88,18) {};
    \node[circle,draw,fill=black,color=black,inner sep=5.83095/4] (v19) at  (118,27) {};
    \node[circle,draw,fill=black,color=black,inner sep=4.12311/4] (v20) at  (84,20) {};
    
    \node[vertex,above of=v4,node distance=10pt,scale=0.7, white,opacity=0] (w4) {};
    \node[vertex,above of=v6,node distance=10pt,scale=0.7, white,opacity=0] (w6) {};
    \node[vertex,above of=v8,node distance=10pt,scale=0.7, white,opacity=0] (w8) {};
    \node[vertex,above of=v11,node distance=10pt,scale=0.7, white,opacity=0] (w11) {};
    \node[vertex,above of=v13,node distance=10pt,scale=0.7, white,opacity=0] (w13) {};
    \node[vertex,above of=v20,node distance=10pt,scale=0.7, white,opacity=0] (w20) {};
    \node[vertex,above of=v1,node distance=10pt,scale=0.7, white,opacity=0] (w1) {};
    \node[vertex,above of=v5,node distance=10pt,scale=0.7, white,opacity=0] (w5) {};
    \node[vertex,above of=v10,node distance=10pt,scale=0.7, white,opacity=0] (w10) {};
    
     \draw [->, red, line width=0.6pt] (w4) -- (v4);
    \draw [->, red, line width=0.6pt] (w6) -- (v6);
    \draw [->, red, line width=0.6pt] (w8) -- (v8);
    \draw [->, red, line width=0.6pt] (w11) -- (v11);
    \draw [->, red, line width=0.6pt] (w13) -- (v13);
    \draw [->, red, line width=0.6pt] (w20) -- (v20);
    \draw [->, red, line width=0.6pt] (w1) -- (v1);
    \draw [->, red, line width=0.6pt] (w5) -- (v5);
    
    \draw [dashed, thick,gray] (v1) -- (v5);
    \draw [dashed, thick,gray] (v1) -- (v7);
    \draw [dashed, thick,gray] (v1) -- (v9);
    \draw [dashed, thick,gray] (v1) -- (v15);
    \draw [dashed, thick,gray] (v2) -- (v7);
    \draw [dashed, thick,gray] (v2) -- (v15);
    \draw [dashed, thick,gray] (v2) -- (v17);
    \draw [dashed, thick,gray] (v3) -- (v7);
    \draw [dashed, thick,gray] (v3) -- (v9);
    \draw [dashed, thick,gray] (v3) -- (v17);
    \draw [dashed, thick,gray] (v5) -- (v7);
    \draw [dashed, thick,gray] (v5) -- (v9);
    \draw [dashed, thick,gray] (v5) -- (v10);
    \draw [dashed, thick,gray] (v5) -- (v15);
    \draw [dashed, thick,gray] (v5) -- (v16);
    \draw [dashed, thick,gray] (v5) -- (v17);
    \draw [dashed, thick,gray] (v7) -- (v9);
    \draw [dashed, thick,gray] (v7) -- (v10);
    \draw [dashed, thick,gray] (v7) -- (v12);
    \draw [dashed, thick,gray] (v7) -- (v14);
    \draw [dashed, thick,gray] (v7) -- (v15);
    \draw [dashed, thick,gray] (v7) -- (v16);
    \draw [dashed, thick,gray] (v7) -- (v17);
    \draw [dashed, thick,gray] (v7) -- (v18);
    \draw [dashed, thick,gray] (v7) -- (v19);
    \draw [dashed, thick,gray] (v9) -- (v10);    
    \draw [dashed, thick,gray] (v9) -- (v12); 
    \draw [dashed, thick,gray] (v9) -- (v14); 
    \draw [dashed, thick,gray] (v9) -- (v15); 
    \draw [dashed, thick,gray] (v9) -- (v16); 
    \draw [dashed, thick,gray] (v9) -- (v17); 
    \draw [dashed, thick,gray] (v9) -- (v18); 
    \draw [dashed, thick,gray] (v9) -- (v19); 
    \draw [dashed, thick,gray] (v10) -- (v12); 
    \draw [dashed, thick,gray] (v10) -- (v14); 
    \draw [dashed, thick,gray] (v10) -- (v15); 
    \draw [dashed, thick,gray] (v10) -- (v17); 
    \draw [dashed, thick,gray] (v10) -- (v18); 
    \draw [dashed, thick,gray] (v10) -- (v19); 
    \draw [dashed, thick,gray] (v12) -- (v15); 
    \draw [dashed, thick,gray] (v12) -- (v15);
    \draw [dashed, thick,gray] (v12) -- (v18);
    \draw [dashed, thick,gray] (v14) -- (v18);
    \draw [dashed, thick,gray] (v15) -- (v16);
    \draw [dashed, thick,gray] (v15) -- (v17);
    \draw [dashed, thick,gray] (v15) -- (v18);
     \draw [dashed, thick,gray] (v17) -- (v18);
     \draw [dashed, thick,gray] (v18) -- (v19);
    \end{tikzpicture}
    }
     
     \vspace{0.2cm}
    \begin{tikzpicture}[scale=0.73][font=\Large]
    \begin{axis}[ybar, enlargelimits=0.23, ymax =20, ymin=3.2,
    ylabel={\textbf{MCN}}, 
    symbolic x coords={A,B,C},
    xtick=data, ylabel near ticks, xticklabels={\textbf{Fed},\textbf{Fast},\textbf{Re-fed}}, title=\textbf{D}, xlabel=\textbf{Phase}, legend pos = north east,legend style={draw=none}, xtick pos=left,ytick pos=left,legend columns=-1,
    ]
\addplot+[blue, blue, fill=blue, opacity=.7] coordinates {(A,12) (B,18) (C,8)};
\addlegendentry{\textbf{Hypergraphs}}

\addplot+[orange, red, fill=red, opacity=.7] coordinates {(A,20) (B,20) (C,20)};
\addlegendentry{\textbf{Graphs}}

\end{axis}
\end{tikzpicture}
\hspace{1cm}
    \caption{Mouse neuron endomicroscopy features. (A), (B) and (C) Neuronal activity networks of the three phases - fed, fast and re-fed, which depicts the spatial location and size of individual cells. Each 2-simplex (i.e., a triangle) represents a hyperedge, and red arrows indicate those control nodes. (D) MCN for the neuronal activity networks modelled by 3-uniform hypergraphs and standard graphs. The cutoff threshold is 0.95 for both the hypergraph and graph models.}
    \label{fig:14}
\end{figure}

The results are shown in Fig. \ref{fig:14}, in which (A), (B) and (C) are network diagrams modelled by 3-uniform hypergraphs for a representative mouse depicting the spatial location and size of individual cells (every 2-simplex is a hyperedge). It is evident from Fig. \ref{fig:14} D that the hypergraph MCN (blue) can successfully differentiate the three phases with different food treatments under the cutoff threshold 0.95. In particular, the fasting phase requires more neurons to be ``controlled'' due to fewer connections, while in the re-fed phase, the MCN is significantly reduced because of an outburst of neuron interactions. On the other hand, the MCN (red in Fig. \ref{fig:14} D) computed from the graph model fails to capture the changes in neuronal activity. The threshold 0.95 in the graph model is too high to produce any connection. This supports the fact that more than two neurons synchronize, or ``co-fire'', in the mouse hypothalamus, and the interactions can be more accurately captured by hypergraphs. Note that our choice of the prescribed threshold is arbitrary, though higher values are desirable as they capture stronger neuronal interactions as the relevant edges/hyeredges. To assess the sensitivity, we also performed our MCN analysis for values of threshold in the range from 0.90 to 0.95 and found similar qualitative results.

\subsection{Allele-Specific Chromosomal Conformation Capture}\label{sec:3.4}

Studies have revealed that there is significant coordination between allelic gene expression biases and local genome architectural changes \cite{lindsly2020functional}. The unbiased genome-wide technology of chromosome conformation capture (Hi-C) has been used to capture the  architecture of the genome through the cell cycle \cite{Rajapakse711,RIED20171,hartwell1989checkpoints}. The notion of transcription factories supports the existence of simultaneous interactions involving multiple genomic loci \cite{cook2018transcription}, implying that the human genome configuration can be represented by a hypergraph \cite{9119161}. In the example, we are given Hi-C data for a small region of chromosome 15 (100 kb bin resolution) which contains two imprinted genes (\textit{SNRPN} and \textit{SNURF}). Imprinted genes are known to only express from one allele, so we want to explore any corresponding differences in local genome architecture \cite{reik2001genomic}. Here we use 4-uniform hypergraphs to partially recover the 3D configuration of the genome according to the multi-correlation (\ref{eq:99}) from the Hi-C matrices.
 
 \begin{figure}[t]
\centering
\tcbox[colback=mygray,top=5pt,left=0pt,right=5pt,bottom=2pt]{
\includegraphics[scale=0.39]{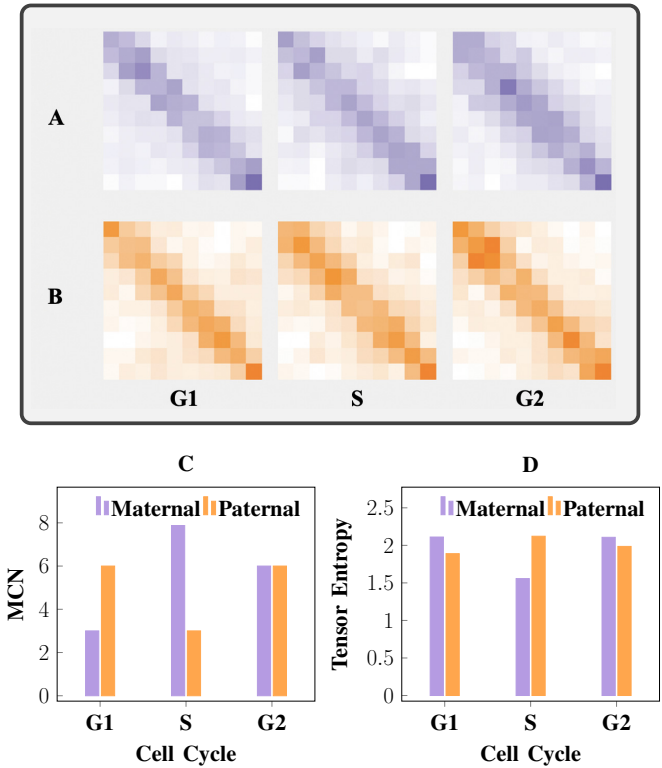}
}

\centering
\vspace{0.2cm}
\begin{tikzpicture}[scale=0.5][font=\LARGE]
    \begin{axis}[ybar, enlargelimits=0.25, ymax =8, ymin=1.4,
    ylabel={\textbf{MCN}}, 
    symbolic x coords={A,B,C},
    xtick=data, ylabel near ticks, xticklabels={\textbf{G1},\textbf{S},\textbf{G2}}, title=\textbf{C}, xlabel=\textbf{Cell Cycle}, legend pos = north east,legend style={draw=none}, xtick pos=left,ytick pos=left,legend columns=-1,
    ]
\addplot+[blue, mypurple, fill=mypurple, opacity=.7] coordinates {(A,3) (B,8) (C,6)};
\addlegendentry{\textbf{Maternal}}

\addplot+[orange, mark options={draw=orange,fill=orange}, opacity=.7] coordinates {(A,6) (B,3) (C,6)};
\addlegendentry{\textbf{Paternal}}

\end{axis}
\end{tikzpicture}
\begin{tikzpicture}[scale=0.5][font=\LARGE]
    \begin{axis}[ybar, enlargelimits=0.25, ymax =2.3, ymin=0.4,
    ylabel={\textbf{Tensor Entropy}}, 
    symbolic x coords={A,B,C},
    xtick=data, ylabel near ticks, xticklabels={\textbf{G1},\textbf{S},\textbf{G2}}, title=\textbf{D}, xlabel=\textbf{Cell Cycle}, legend pos = north east,legend style={draw=none}, xtick pos=left,ytick pos=left,legend columns=-1,
    ]
\addplot+[blue, mypurple, fill=mypurple, opacity=.7] coordinates {(A,2.1109) (B,1.5572) (C,2.1072)};
\addlegendentry{\textbf{Maternal}}

\addplot+[orange, mark options={draw=orange,fill=orange}, opacity=.7] coordinates {(A,1.8893) (B,2.1216) (C,1.9860)};
\addlegendentry{\textbf{Paternal}}

\end{axis}
\end{tikzpicture}

\caption{Allele-specific Hi-C features. (A) and (B) Hi-C maps of a local region surrounding the imprinted genes \textit{SNRPN} and \textit{SNURF} from the maternal and paternal Chromosome 15, respectively, through the cell cycle phases G1, S and G2. The darker the color, the more interactions between two loci. (C) MCN of the 4-uniform hypergraphs, recovered from Hi-C measurements with multi-correlation cutoff threshold 0.99, through the cell cycle phases G1, S and G2. (D) Tensor entropies of the 4-uniform hypergraphs described in (C).}
\label{fig:15}
\end{figure}
 
The results are shown in Fig. \ref{fig:15}. Clearly, it is hard to tell the difference between the maternal and paternal genome architectures directly from the Hi-C maps, see Fig. \ref{fig:15} A and B. However, after converting to hypergraphs, we can easily detect the structural discrepancy using the notion of MCN in the cell cycle phases G1 and S, see Fig. \ref{fig:15} C. Although the MCN are equal between the maternal and paternal networks in the cell cycle phase G2, the maximum possible choices of minimum subsets of control nodes are different (one is twenty four, and one is thirty three). This indicates that there are some architectural similarities between the maternal and paternal architectures in G2 compared to the previous two phases. Furthermore, we corroborate our results by using the notion of tensor entropy. Tensor entropy is a spectral measure, which can decipher topological attributes of uniform hypergraphs \cite{9119161}. In particular, the two results of tensor entropy and MCN are consistent, in the sense that the largest gap of tensor entropy between the maternal and paternal architectures occurs in S, and the smallest gap occurs in G2, see Fig. \ref{fig:15} D. Biologically, in S phase, DNA replication of the genomes may lead to a large structural dissimilarity between the maternal and paternal architectures, while in the G2 phase, both genomes prepare for mitosis which may result in a small structural discrepancy. Moreover, we believe that the control loci (nodes) can play a significant role in cellular reprogramming, a process that introduces proteins called transcription factors as a control mechanism for transforming one cell type into another. 
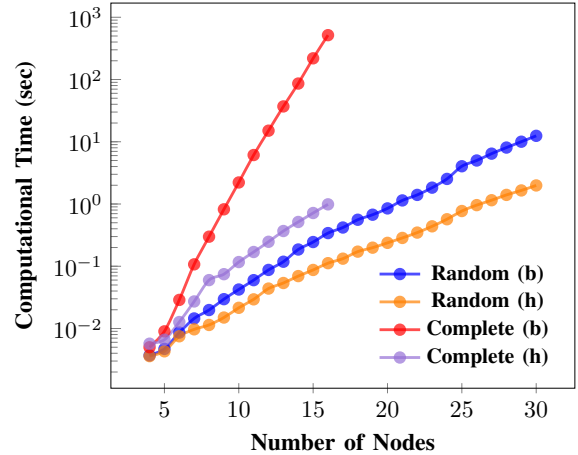
\begin{figure}[t!]
\centering
\begin{tikzpicture}[scale=0.9]
\begin{axis}[ymode=log, ytick={1e-5, 1e-4, 1e-3, 1e-2, 1e-1, 1, 10, 100,1000}, xlabel={\textbf{Number of Nodes}}, ylabel={\textbf{Computational Time (sec)}}, ylabel near ticks, legend pos=south east, xlabel near ticks,xtick pos=left,ytick pos=left,xtick align=center,legend style={draw=none, font=\small} ]

\addplot+[smooth, very thick,mark=*,blue, mark options={draw=blue,fill=blue}, opacity=.7] coordinates {(4,0.0037) (5,0.0047) (6,0.0086) (7,0.0145) (8,0.0198) (9,0.0296) (10,0.0423) (11,0.0598) (12,0.0873) (13,0.1182) (14,0.1858) (15,0.2449) (16, 0.3410) (17, 0.4190) (18, 0.5586) (19, 0.6730) (20, 0.8497) (21, 1.1419) (22, 1.4003) (23, 1.8310) (24, 2.5292) (25, 4.0413) (26, 4.9987) (27, 6.4402) (28,8.0666) (29,10.0164) (30,12.4514)};
\addlegendentry{\textbf{Random (b)}}

\addplot+[smooth, very thick,mark=*,orange, mark options={draw=orange,fill=orange}, opacity=.7] coordinates {(4,0.0036) (5,0.0043) (6,0.0075) (7,0.0097) (8,0.0114) (9,0.0150) (10,0.0215) (11,0.0294) (12,0.0438) (13,0.0538) (14,0.0695) (15,0.0872) (16, 0.1121) (17, 0.1324) (18, 0.1717) (19,0.1989) (20, 0.2376) (21, 0.2850) (22, 0.3473) (23, 0.4380) (24, 0.5691) (25, 0.7715) (26,0.9595) (27,1.1503) (28,1.4005) (29,1.6520) (30,1.9871)};
\addlegendentry{\textbf{Random (h)}}

\addplot+[smooth, very thick,mark=*,red, mark options={draw=red,fill=red}, opacity=.7] coordinates {(4,0.0050) (5,0.0090) (6,0.0288) (7,0.1073) (8,0.2985) (9,0.8217) (10,2.2188) (11,6.1192) (12,15.0671) (13,36.9600) (14,86.1357) (15,219.0578) (16,515.6578)};
\addlegendentry{\textbf{Complete (b)}}

\addplot+[smooth, very thick,mark=*,mypurple, mark options={draw=mypurple,fill=mypurple}, opacity=.7] coordinates {(4,0.0057) (5,0.0066) (6,0.0127) (7,0.0272) (8,0.0599) (9,0.0746) (10,0.1168) (11,0.1694) (12,0.2472) (13,0.3682) (14,0.5146) (15,0.7168) (16,0.9876)};
\addlegendentry{\textbf{Complete (h)}}

\end{axis}
\end{tikzpicture}
\hspace{0.9cm}
\caption{Computational time comparisons in determining MCN. In the legend, the letter b stands for the brute-force search, while the letter h stands for the heuristic approach. Since the computational time using a brute-force search in determining the MCN of the complete uniform hypergraphs grows very fast, we only compute them up to sixteen nodes for comparison. For the purpose of accuracy, we ran each algorithm five times and took the average of the computational times.}
\label{fig:23}
\end{figure}

\subsection{MCN Computation Comparison}\label{sec:3.5}
In this example, we compare the computational efficiency of the heuristic approach described in Algorithm \ref{alg:2} and brute-force search in finding the MCN of random 4-uniform hypergraphs (with hyperedge density 50\%) and complete 4-uniform hypergraphs. The results are shown in Fig. \ref{fig:23}. Evidently, Algorithm \ref{alg:2} is more time-efficient than a brute-force search as the number of nodes becomes larger in both the configurations. In particular, when a uniform hypergraph is complete (or nearly complete), the heuristic exhibits a huge time advantage. More significantly, it produces exactly the same MCN as a brute-force search in these simulations.

\section{Discussion}\label{sec:4}
The first four numerical examples reported here highlight that the tensor-based multilinear system (\ref{eq:2}) can characterize the multidimensional interactions in uniform hypergraphs. The MCN of uniform hyperchains, hyperrings and hyperstars, and their variants are related to their degree distributions. It is also a good indicator of uniform hypergraph robustness. However, more theoretical and numerical investigations are required to evaluate the controllability of more general uniform hypergraphs, and its relation to the hypergraph topology. Moreover, in practice, hypergraphs like co-authorship networks and protein-protein interaction networks are very large, so computing the reduced controllability matrix and its corresponding MCN is still challenging. Tensor decomposition and a ``maximum matching'’ type approach need to be considered in order to facilitate efficient computations \cite{doi:10.1137/07070111X,hopcroft1973n,BaKo06}.

Instead of looking at uniform hypergraphs, can we think of controllability of more general hypergraphs where each hyperedge contains an arbitrary amount of nodes? The main idea is to make non-uniform hypergraphs uniform, which can then be represented by tensors. In the following, we adopt the definition of generalized adjacency tensors of non-uniform hypergraphs from \cite{BANERJEE201714}.

\textit{Definition 9 (\cite{BANERJEE201714}):} Let \textsf{G} = \{\textbf{V}, \textbf{E}\} be a hypergraph with $n$ nodes, and $k$ be the maximum cardinality of the hyperedges. The adjacency tensor $\textsf{A}\in\mathbb{R}^{n\times n\times \dots\times n}$ of \textsf{G}, which is a $k$-th order $n$-dimensional supersymmetric tensor, is defined as
\begin{equation}
\textsf{A}_{j_1j_2\dots j_k} = \begin{cases} \frac{s}{\alpha} \text{ if $(i_1,i_2,\dots,i_s)\in \textbf{E}$}\\ \\0, \text{ otherwise}\end{cases},
\end{equation}
where, $j_l\in\{i_1,i_2,\dots,i_s\}$ with at least once for each element of the set for $l=1,2,\dots,k$, and 
\begin{equation*}
    \alpha =\sum_{k_1+k_2+\dots+k_s=k} \frac{k!}{\prod_{l=1}^s k_i!}.
\end{equation*}

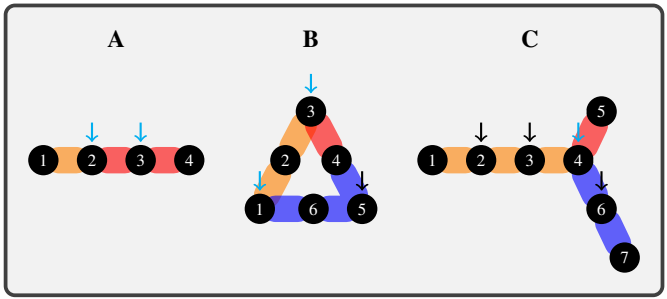
\begin{figure}[t]
\centering
\tcbox[colback=mygray,top=0pt,left=5pt,right=5pt,bottom=5pt]{
\begin{tikzpicture}[scale=0.92, transform shape]
\node[vertex,text=white,scale=0.7] (v1) {1};
\node[vertex,right of=v1,node distance = 20pt,text=white,scale=0.7] (v2) {2};
\node[vertex,right of=v2,node distance = 20pt,text=white,scale=0.7] (v3) {3};
\node[vertex,right of=v3,node distance = 20pt,text=white,scale=0.7] (v4) {4};

\node[vertex,above of=v2,node distance=20pt,mygray,scale=0.7] (w1) {};
\node[vertex,above of=v3,node distance=20pt,mygray,scale=0.7] (w2) {};
\path [->,shorten >=1pt,shorten <=1pt, thick, cyan](w1) edge node[left] {} (v2);
\path [->,shorten >=1pt,shorten <=1pt, thick,cyan](w2) edge node[left] {} (v3);

\node[vertex,below right = 11pt and 20pt of v4,text=white,scale=0.7,node distance=20pt] (v5) {1};
\node[vertex,above right = 11pt and 1.5pt of v5,text=white,scale=0.7] (v6) {2};
\node[vertex,above right = 11pt and 1.5pt of v6,text=white,scale=0.7] (v7) {3};
\node[vertex,below right = 11pt and 1.5pt of v7,text=white,scale=0.7] (v8) {4};
\node[vertex,below right = 11pt and 1.5pt of v8,text=white,scale=0.7] (v9) {5};
\node[vertex,left of=v9,text=white,scale=0.7, node distance = 20pt] (v10) {6};

\node[vertex,above of=v7, node distance = 30pt,mygray, text=black,scale=1] (b) {\textbf{B}};
\node[vertex,left of=b, node distance = 80pt,mygray, text=black,scale=1] (a) {\textbf{A}};
\node[vertex,right of=b, node distance = 90pt,mygray, text=black,scale=1] (c) {\textbf{C}};

\node[vertex,above of=v5,node distance=20pt,mygray,scale=0.7] (w3) {};
\node[vertex,above of=v7,node distance=20pt,mygray,scale=0.7] (w4) {};
\node[vertex,above of=v9,node distance=20pt,mygray,scale=0.7] (w5) {};
\path [->,shorten >=1pt,shorten <=1pt, thick,cyan](w3) edge node[left] {} (v5);
\path [->,shorten >=1pt,shorten <=1pt, thick,cyan](w4) edge node[left] {} (v7);
\path [->,shorten >=1pt,shorten <=1pt, thick](w5) edge node[left] {} (v9);

\node[vertex,right of=v4,text=white,scale=0.7,node distance=120pt] (v11) {2};
\node[vertex,left of=v11,text=white,scale=0.7,node distance=20pt] (v15) {1};
\node[vertex,right of=v11,text=white,scale=0.7,node distance = 20pt] (v12) {3};
\node[vertex,right of=v12,text=white,scale=0.7,node distance = 20pt] (v13) {4};
\node[vertex,above right = 11pt and 0.5pt of v13,text=white,scale=0.7] (v14) {5};
\node[vertex,below right = 11pt and 0.5pt of v13,text=white,scale=0.7] (v16) {6};
\node[vertex,below right = 11pt and 0.5pt of v16,text=white,scale=0.7] (v17) {7};

\node[vertex,above of=v11,node distance=20pt,mygray,scale=0.7] (w6) {};
\node[vertex,above of=v12,node distance=20pt,mygray,scale=0.7] (w7) {};
\node[vertex,above of=v13,node distance=20pt,mygray,scale=0.7] (w8) {};
\node[vertex,above of=v16,node distance=20pt,mygray,scale=0.7] (w9) {};
\path [->,shorten >=1pt,shorten <=1pt, thick](w6) edge node[left] {} (v11);
\path [->,shorten >=1pt,shorten <=1pt, thick](w7) edge node[left] {} (v12);
\path [->,shorten >=1pt,shorten <=1pt, thick,cyan](w8) edge node[left] {} (v13);
\path [->,shorten >=1pt,shorten <=1pt, thick](w9) edge node[left] {} (v16);

\begin{pgfonlayer}{background}
\draw[edge,color=orange] (v1) -- (v2);
\draw[edge,color=red] (v2) -- (v3) -- (v4);

\draw[edge,color=orange] (v5) -- (v6) -- (v7);
\draw[edge,color=red] (v7) -- (v8);
\draw[edge,color=blue] (v8) -- (v9) -- (v10) -- (v5);

\draw[edge,color=orange] (v15) -- (v11) -- (v12) -- (v13);
\draw[edge,color=red] (v13) -- (v14);
\draw[edge,color=blue] (v13) -- (v16) -- (v17);

\end{pgfonlayer}
\end{tikzpicture}
}
\caption{Controllability of non-uniform hypergraphs with MCN. (A) Non-uniform hyperchain with $e_1=\{1,2\}$ and $e_2=\{2,3,4\}$. (B) Non-uniform hyperring with $e_1=\{1,2,3\}$, $e_2=\{3,4\}$ and $e_3=\{4, 5, 6, 1\}$. (C) Non-uniform hyperstar with $e_1=\{1,2,3,4\}$, $e_2=\{4,5\}$ and $e_3=\{4,6,7\}$. Nodes with arrows from the top are the control nodes, and the cyan arrows indicate the control nodes with the highest degrees in the configurations.
}
\label{fig:20}
\end{figure}

The choice of the nonzero coefficients $\frac{s}{\alpha}$ preserves the degree of each node, i.e., the degree of node $j$ computed using  (\ref{eq:r6}) with weights as defined above is equal to number of hyperedges containing the node in the original non-uniform hypergraph. When \textsf{G} is uniform, the above definition reduces to Definition 1. See \cite{BANERJEE201714} for examples. Once we have the adjacency tensor of a hypergraph, we can discuss the controllability of the hypergraph when $k$ is even using the techniques developed in Section \ref{sec:2.4}. For curiosity,  we build several non-uniform hypergraphs and determine their MCN. The results are shown in Fig. \ref{fig:20}. Intriguingly, the control strategies for non-uniform hyperchains, hyperrings and hyperstars with one overlapping nodes between hyperedges are similar to those discussed in Section \ref{sec:3.2}. High degree nodes are preferred to be controlled with each hyperedge containing $s-1$ control nodes where $s$ is the cardinality of the hyperedge (there is one exception in the hyperring configuration). 

Furthermore, one may consider weights and directions for hypergraphs. Weighted hypergraphs can be readily accomplished by replacing $\frac{1}{(k-1)!}$ with $\textsf{W}_{j_1j_2\dots j_k}$ for some supersymmetric weight tensors \textsf{W} in (\ref{eq:98}). Banerjee \textit{et al.} \cite{banerjee2017spectrum} also defined adjacency tensors for \textit{directed hypergraphs}, so it will be interesting to establish the controllability of directed hypergraphs.

\section{Conclusion}\label{sec:5}
In this paper, we proposed a new notion of controllability for uniform hypergraphs based on tensor algebra and polynomial control theory. We represented the dynamics of uniform hypergraphs by a tensor product based multilinear system, and derived a Kalman-rank-like condition to determine the controllability of even uniform hypergraphs. We established theoretical results on the MCN of even uniform hyperchains, hyperrings and hyperstars as well as complete even uniform hypergraphs. We proposed MCN as a measure of hypergraph robustness, and found that is it related to the hypergraph degree distribution. We also presented a heuristic to obtain the MCN efficiently. Additionally, we applied the notion of MCN to the real world biological networks to quantify structural differences, and achieved outstanding performances. Finally, we discussed the controllability of general hypergraphs. As mentioned in Section \ref{sec:4}, a ``maximum matching'' type approach is needed in order to find the MCN  of large undirected/directed hypergraphs in a scalable fashion. On the other hand, more work is required to fully understand the control properties of the tensor-based multilinear system (\ref{eq:2}). For example, it will be useful to realize the potential of tensor algebra based computations for controllability Gramians and Lyapunov equations.  In addition, it will be worthwhile to develop theoretical and computational frameworks for observer and feedback control design, and apply them to the dynamics of hypergraphs. 

\section*{Acknowledgments}
We would like to thank Dr. Frederick Leve at the Air Force Office of Scientific Research (AFOSR) for support and encouragement. We would also like to thank the three referees for their constructive comments, which lead to a significant improvement of the paper.

\ifCLASSOPTIONcaptionsoff
  \newpage
\fi



\bibliographystyle{IEEEtran}
\bibliography{tensorref}
%


%

\begin{IEEEbiography}[{\includegraphics[width=1in,height=1.25in,clip,keepaspectratio]{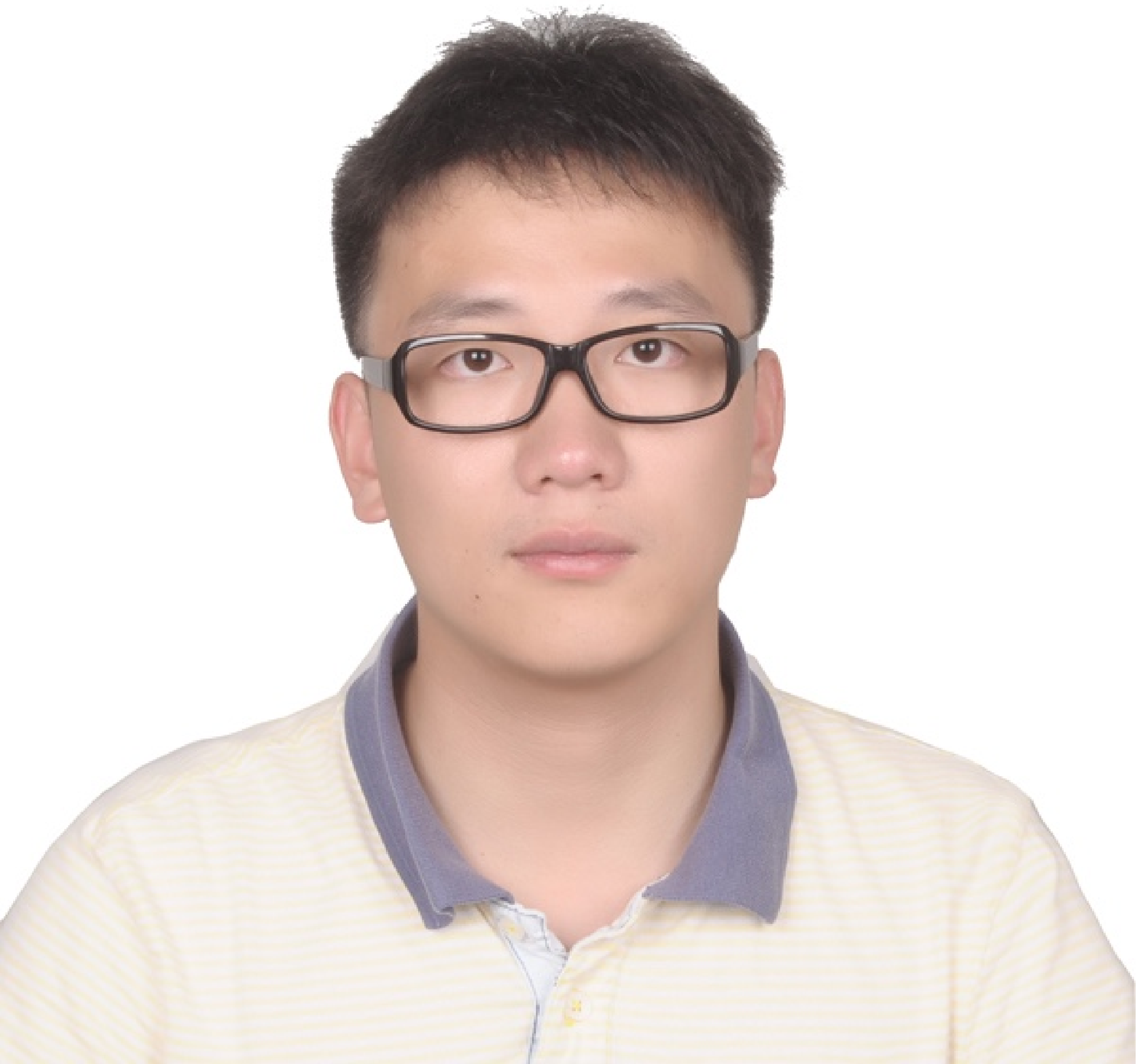}}]{Can Chen} is currently a Ph.D. candidate in Applied \& Interdisciplinary Mathematics and a Master's student in Electrical \& Computer Engineering at the University of Michigan, Ann Arbor. He received B.S. degree in Mathematics at the University of California, Irvine in 2016. His research is focused on data-guided control of multiway dynamical systems. 
\end{IEEEbiography}

\begin{IEEEbiography}[{\includegraphics[width=1in,height=1.25in,clip,keepaspectratio]{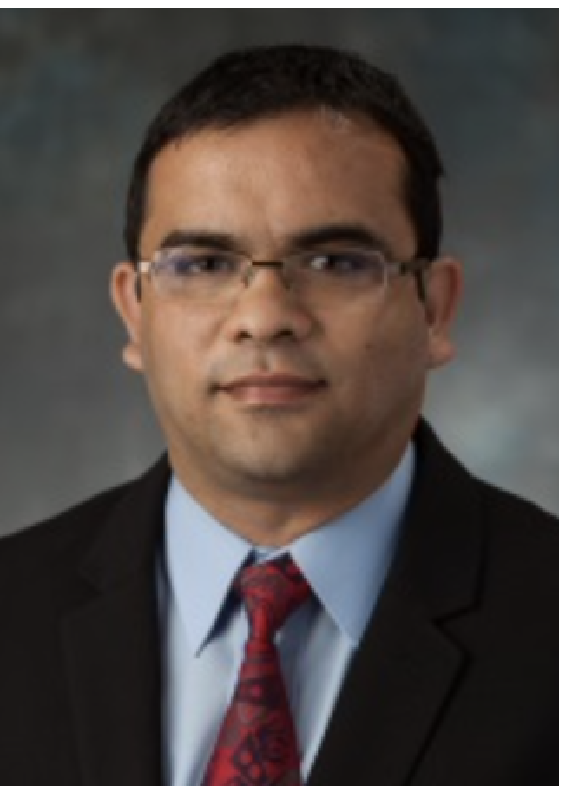}}]{Amit Surana} is currently a technical fellow at Raytheon Technologies Research Center. He received his Bachelor’s degree in Mechanical Engineering from Indian Institute of Technology Bombay in 2000, his M.S. in Mechanical Engineering and M.A. in Mathematics both from Pennsylvania State University in 2002 and 2003, respectively, and his Ph.D. in Mechanical Engineering from Massachusetts Institute of Technology (MIT) in 2007.  His research interests include dynamical systems, control theory, machine learning, and multi-agent systems with a broad range of applications in aerospace and defense domains.
\end{IEEEbiography}

\begin{IEEEbiography}[{\includegraphics[width=1in,height=1.25in,clip,keepaspectratio]{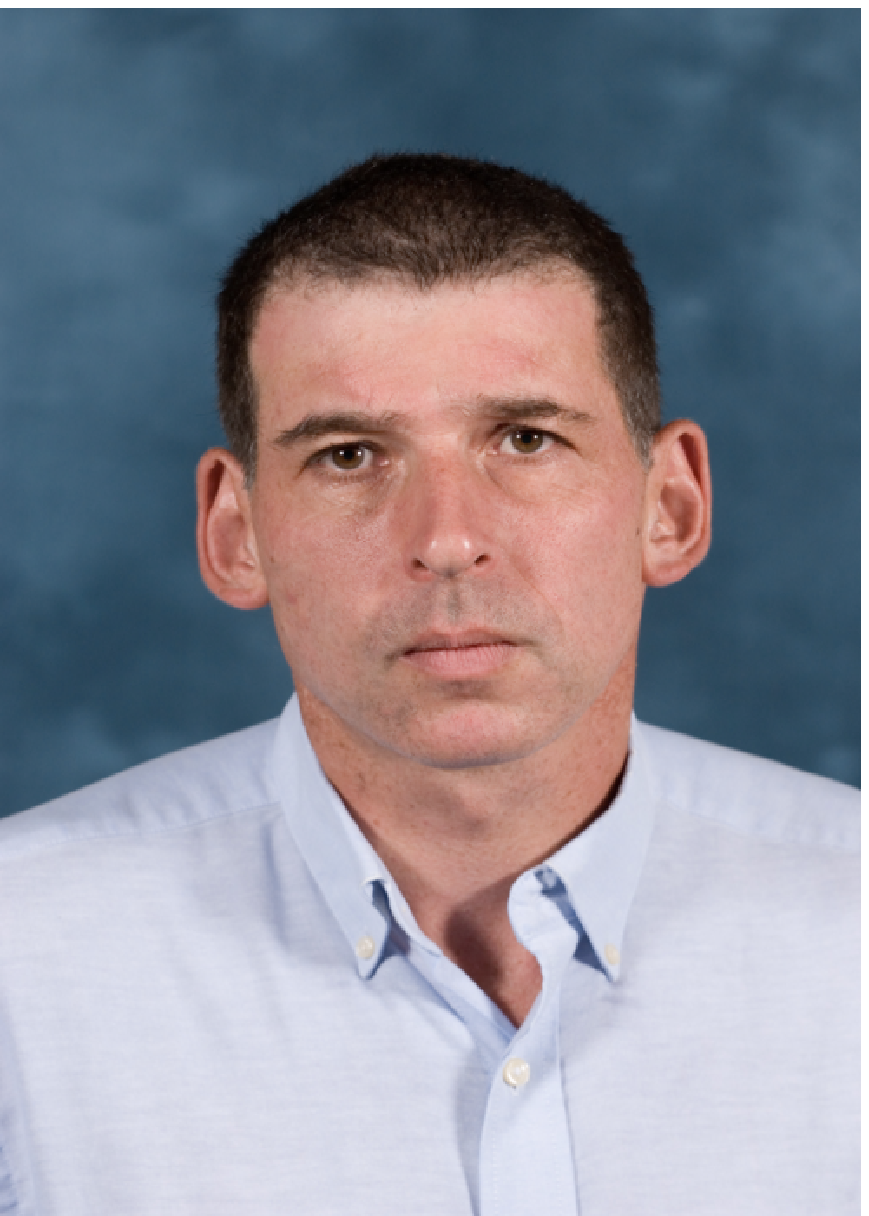}}]{Anthony M. Bloch} is the Alexander Ziwet Collegiate
Professor of Mathematics at the University of Michigan
where he is currently department chair. 
He received a B.Sc. (Hons) in Applied Mathematics and Physics
from the University of the Witwatersrand, Johannesburg, in 1978, an M.S.
in Physics from the California Institute of Technology in 1979, an M. Phil
in Control Theory and Operations Research from Cambridge University in
1981 and a Ph.D. in Applied Mathematics from Harvard University in 1985. He is Fellow 
of the IEEE, SIAM and the AMS. He has served on the editorial boards of various journals, was Editor-in-Chief of the SIAM Journal 
of Control and Optimization and is currently coEditor-in-Chief of the Journal
of Nonlinear Science. 
\end{IEEEbiography}

\begin{IEEEbiography}[{\includegraphics[width=1in,height=1.25in,clip,keepaspectratio]{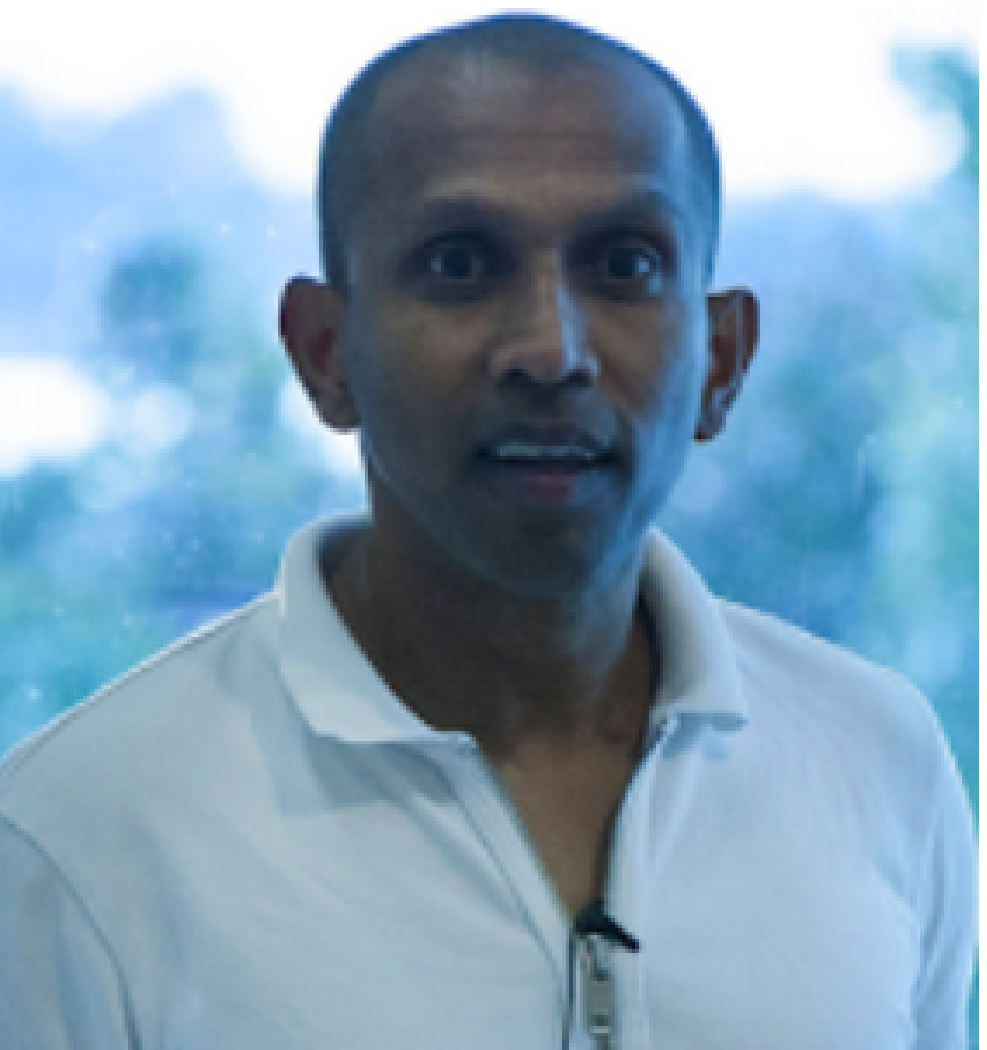}}]{Indika Rajapakse} is currently an Associate Professor of Computational Medicine \& Bioinformatics, in the Medical School, and an Associate Professor of Mathematics at the University of Michigan, Ann Arbor. He is also a member of the Smale Institute. His research is at the interface of biology, engineering and mathematics. His areas include dynamical systems, networks, mathematics of data and cellular reprogramming.
\end{IEEEbiography}
\vfill




\end{document}